\documentclass{article}

\usepackage[english]{babel}
\usepackage[a4paper,top=3cm,bottom=2cm,left=2cm,right=2cm,marginparwidth=1.75cm]{geometry}

\usepackage[style = numeric, sorting = nty]{biblatex}

\usepackage{float}
\usepackage{wrapfig}
\usepackage{amsmath}
\usepackage{xcolor}
\usepackage{mathtools}
\usepackage{authblk}
\usepackage{enumitem}
\newtheorem{remark}{Remark}
\newcommand{\rv}[1]{\textcolor{black}{#1}}
\newcommand{\vr}[1]{\textcolor{black}{#1}}
\usepackage{csquotes}

\newif\ifshowfig 
 \showfigtrue      

\title{Numerical solution of two dimensional scalar conservation laws using compact implicit numerical schemes on Cartesian meshes \thanks{The work was supported by the grants VV-MVP-24-0116, VEGA 1/0314/23, and APVV-23-016.}}

\author[1]{Peter FROLKOVI\v{C}}
\affil[1]{Department of Mathematics and Descriptive Geometry, Slovak University of Technology in Bratislava (peter.frolkovic@stuba.sk)}
\author[2]{Dagmar \v{Z}\'AKOV\'A}
\affil[2]{Department of Mathematics and Descriptive Geometry, Slovak University of Technology in Bratislava (dagmar.zakova@stuba.sk)}

\begin{document}

\date{}
\maketitle

\textbf{Abstract: }  
This paper deals with the numerical solution of conservation laws in the two dimensional case using a novel compact implicit time discretization that enables applications of fast algebraic solvers. We present details for the second order accurate parametric scheme based on the finite volume method including simple variants of ENO (Essentially Non-Oscillatory) and WENO (Weighted Essentially Non-Oscillatory) approximations. \rv{To avoid oscillatory numerical solutions for large time steps, we propose limiting in time which is provable non-oscillatory in the case of linear advection equation with variable velocity.} We present numerical experiments for representative linear and nonlinear problems.

\textbf{Keywords: }{compact implicit discretization, high-resolution numerical method, positive coefficients scheme, advection equation, Burgers equation}

\section{Introduction}

Time dependent partial differential equations are typically solved by numerical methods where the time discretization is applied separately from the spatial discretization. Such an approach has the advantage that well developed and implemented numerical integration methods for the solutions of ordinary differential equations (ODEs) can be used after the equation is discretized in space. Already at the level of ODEs, it is recognized that the explicit type of time discretization methods is not suitable for systems of ODEs of which the solution describes processes with different time scales. In such a case, the restriction on the size of the time step of explicit schemes is no longer applicable if the fast processes require too small time steps only due to the stability constraints and not due to the accuracy requirements. A similar situation can occur when the time dependent solution in numerical simulations is approaching a stationary state, especially if some asymptotic preserving properties are expected from the numerical solution.

In the above mentioned cases, the well-known remedy is to apply implicit methods for the numerical solution of the ODEs that replace some explicit formulas to determine the numerical solution by discrete algebraic equations to be solved. This type of method is available with enhanced stability properties and offers, in an ideal case, the usage with the time steps having a size that is limited only by the accuracy requirements. The price to pay for such stability enhancement is strongly dependent on the efficiency with which the resulting discrete systems of algebraic equations can be solved. 

Numerical methods based on the separate discretization of time and space have been well developed for hyperbolic systems.
Concerning implicit methods (or their combinations with explicit methods), we can mention the class of Runge-Kutta (RK) methods used for hyperbolic problems in \cite{pareschi2005implicit, boscarino_high_2022,gomez-buenoImplicitSemiimplicitWellbalanced2023,condeImplicitImplicitExplicit2017, arbogast2020third, puppo_quinpi_2022, michel2022tvd, birkenConservationPropertiesIterative2022,puppo2023quinpi}. Furthermore, it is recognized that further improvements can be achieved if, opposite to standard RK methods, not only the values of the right hand side function are used, but also its derivatives. This approach is realized in the development of two- or multi-derivative RK methods \cite{seal_high-order_2014,tsai_two-derivative_2014} available also in the form of implicit methods \cite{zeifang_two-derivative_2021,gottlieb_high_2022}. Although all these methods are very well developed, their applications are not straightforward if some structural properties of the solution for hyperbolic equations shall be preserved not only by space discretization, but also by the time discretization.

The main idea of multi-derivative RK methods is shared with other types of time discretization methods that exploit finite Taylor series expansions of the solution with the so-called Lax-Wendroff (LW) or Cauchy-Kowalevskaya procedure when the time derivatives are replaced by corresponding derivatives of the right hand side function. 
The methods based on the LW procedure offer a great opportunity to couple the time and space discretizations, so these methods can be developed in a hybrid way in which it can be easier to obtain some properties of numerical schemes.
The simplest tool to realize it is to approximate each term in the Taylor series using different space discretizations \cite{qiu_finite_2003, seal_high-order_2014, tsai_two-derivative_2014, li_two-stage_2019,zeifang_two-derivative_2021}. The coupling of time and space discretization is strengthened in methods that use mixed spatial-temporal derivatives in the Lax-Wendroff procedure \cite{zorio_approximate_2017, frolkovic2018semi, carrillo2021lax,maccaSemiImplicittypeOrderadaptiveCAT22023}.

As shown in \cite{frolkovic2018semi,frolkovic_semi-implicit_2022,frolkovic2023high}, the idea of implicit time discretization using the LW procedure involving the mixed spatial-temporal derivatives can be used not only to enhance the stability of the method but also to improve the solvability of the resulting algebraic systems. In particular, a compact high-resolution finite difference method for hyperbolic problems in the one-dimensional case is presented in \cite{frolkovic2023high} in connection with the fractional step method \cite{lozano_implicit_2021} where the Jacobian matrix has a convenient structure, and therefore efficient solvers like the fast sweeping method \cite{lozano_implicit_2021} can be used. 

\rv{A unique feature to implicit time discretizations is so called limiting in time to avoid oscillatory numerical solutions \cite{duraisamy_implicit_2007,arbogast2020third,puppo_quinpi_2022,michel2022tvd,puppo2023quinpi} as it is typically necessary only for large time steps that are not allowed anyway for explicit schemes due to stability reasons. Clearly, discontinuities in the solutions of (nonlinear) hyperbolic problems can occur not only along space variable, but also in time, therefore, analogously to well developed limiters in space \cite{harten_class_1984,sweby1984high, leveque_finite_2004} or (W)ENO space reconstructions \cite{shu_essentially_1998}, similar procedure is necessary for the time variable. Unfortunately, the stencils for the approximations in time are not so flexible as in the space, therefore, to avoid oscillatory approximations for discontinuities or rapid changes in time, one typically reduce locally high order approximations to the first order accurate form according to some indicators \cite{arbogast2020third,puppo_quinpi_2022,michel2022tvd,puppo2023quinpi}. It was already noticed in \cite{duraisamy_implicit_2007} that such approach can result in numerical solutions that are corrupted by artificial diffusion. The authors in \cite{duraisamy_implicit_2007} proposed a time limiter that 
depends on so called Courant numbers and that acts eventually only when this numbers are locally larger than one (the typical stability limit of explicit schemes). The price for this improvements is that fully coupled nonlinear algebraic system has to be solved even for the linear advection problems.}

\rv{
In \cite{frolkovic2023high}, the time limiter of \cite{duraisamy_implicit_2007} was modified to obtain the compact implicit high-resolution numerical scheme. Moreover, when used in  \cite{frolkovic2023high} with the fractional time step method \cite{lozano_implicit_2021}, a convenient structure of Jacobian matrix is obtained in each step, and the effective algebraic solvers like the fast sweeping method can be applied.
}

In this work, we extend the approach of \cite{frolkovic2023high} in several aspects which indicate important steps towards applying the compact implicit scheme in even more general settings. The numerical scheme here is based on the two-dimensional finite volume method that offers more flexibility when, e.g., an extension for nonuniform grids is desired \cite{shu_essentially_1998}, while the 1D finite difference method in \cite{frolkovic2023high} is tied up to uniform grids. 
In contrast to \cite{lozano_implicit_2021,frolkovic2023high}, we do not use the fractional time step method, which can in many cases decrease accuracy of results due to time splitting errors \cite{leveque_finite_2004}.  To avoid unphysical oscillations, we propose a high-resolution form of the compact implicit scheme using a Weighted Essentially Non-Oscillatory (WENO) reconstruction in space \cite{shu_essentially_1998} not given in \cite{frolkovic2023high}, and a nonlinear limiting in time using a flux limiter \rv{that is tailored to the fast sweeping method}. Opposite to the 1D case in \cite{frolkovic2023high}, we treat here a space dependent velocity in 2D situations and relate the computations of  parameters \rv{in the time limiter of} the numerical scheme with theoretical results on a preservation of non-negative coefficients in the final numerical scheme.  \rv{Similarly to \cite{duraisamy_implicit_2007,frolkovic2023high}, we adopt this approach also to nonlinear cases, and, additionally, we couple it with the fast sweeping iterations that simplify the nonlinearities due to the solution dependent numerical parameters, so a better behavior of required nonlinear solvers can be observed.}

The main purpose of this study, similarly to \cite{arbogast2020third,puppo_quinpi_2022,frolkovic2023high}, is to show that, following these guidelines, the resulting compact implicit scheme can be used for accurate numerical solutions of representative equations of scalar hyperbolic equations. The main motivation is to show that if the method is used as an implicit ``black-box solver'',  a good accuracy can be obtained in standard numerical experiments that is significantly improved with respect to the implicit first order accurate scheme. Consequently, the dominant criterion for the choice of time steps is the accuracy of numerical resolution for the solution phenomena that are of
interest, and no restriction on the time step is required due to the stability requirement that might be related to some uninteresting phenomena. In this initial study, we do not compute numerical examples in which such stability requirements are presented, which we plan to do in our subsequent research. 

The paper is organized as follows. In Section \ref{sec-model} we introduce the scalar conservation laws for two dimensional case and a general form of the finite volume method for its numerical solution. \rv{In Section \ref{sec-fsm} we introduce the fast sweeping method to solve the resulting algebraic systems together with a short overview of the related literature.} In Section \ref{sec-compact} we derive a parametric form of the second order accurate scheme. \rv{In Section \ref{sec-hr} we introduce the high-resolution form using the ENO (Essentially Non-Oscillatory) and WENO  approximations and time limiters.  In Section 6 we derive the conditions when the high-resolution form produces positive coefficients scheme and define the time limiter accordingly.} In Section \ref{sec-numerical} we illustrate the properties of the presented schemes with several numerical experiments, and we conclude in Section \ref{sec-conclude}.

\section{Mathematical model and finite volume method}
\label{sec-model}

We consider two representative scalar hyperbolic partial differential equations in the two dimensional case. 
First, we consider the nonlinear hyperbolic equation in the form 
\begin{equation}
\label{2D_general1}
    \partial_{t} u +  \nabla \cdot \boldsymbol{f}(u)= 0,
\end{equation}
with $\boldsymbol{f}(u) = (f(u), g(u))$ being the vector flux function, so \eqref{2D_general1} can be written in the form
\begin{equation}
\label{2D_general}
\partial_{t} u +  \partial_{x} f(u) +  \partial_{y} g(u) = 0 
\end{equation}
with $u = u(x,y,t)$ being the unknown function for $t\in(0,T)$ and $x\in(x_L,x_R)\subset R$, $y\in(y_L,y_R)\subset R$. The initial condition is defined by $u(x,y,0)=u^0(x,y)$ and the Dirichlet boundary conditions, if prescribed, are denoted by
\begin{equation}
    \label{boundary2D}
    u(x_L,y,t) = u_{x_L} (y,t), \quad u(x,y_L,t) = u_{y_L} (x,t),  \quad u(x_R,y,t) = u_{x_R} (y,t), \quad u(x,y_R,t) = u_{y_R} (x,t) \,.
\end{equation}

Next, we introduce the numerical discretization to solve the equation \eqref{2D_general}.
The discretization is done in space and time using the following notation. We denote $t^n=n\tau$, $n=0,1,...N$ for a chosen $N$ and $\tau>0$ with $t^{n+1/2}=t^n+\tau/2$. The spatial discretization is based on the finite volume method \cite{eymard2000finite,leveque_finite_2004,barth2018finite}. For simplicity of the notation, we assume a squared computational domain that is divided into finite volumes of the form $I_{i,j} = (x_{i-1/2},x_{i+1/2}) \times (y_{j-1/2},y_{j+1/2})$ with the regular square grid, where $x_{i-1/2} = x_L + (i-1)h$ and $y_{j-1/2} = x_L + (j-1)h$ with $i,j=1, 2, ..., M$ for the chosen $M$ and $h=(x_R-x_L)/M$.

The main idea behind the finite volume method is to integrate the differential equation (\ref{2D_general}) over $I_{i,j} \times (t^n,t^{n+1})$
\begin{equation}
    \label{int2Dgen}
    \int_{t^n}^{t^{n+1}} \int_{I_{i,j}} \big( \partial_{t} u +  \partial_{x} f(u) + \partial_y g(u) \big) \,dx \,dy\,dt = 0.
\end{equation}
One can write \eqref{int2Dgen} formally in the form,
\begin{equation}
\label{exactscheme}
\tilde u_{i,j}^{n+1} -  \tilde u_{i,j}^n + \frac{\tau}{h} ( f_{i+1/2,j} - f_{i-1/2,j}) + \frac{\tau}{h} ( g_{i,j+1/2} - g_{i,j-1/2}) = 0.
\end{equation}
where $\tilde u_{i,j}^n$ denotes the space averaged values and $f_{i+1/2,j}$ and $g_{i,j+1/2}$ the time averaged values over corresponding integrals.

We aim to obtain a second order accurate approximation of \eqref{exactscheme}, therefore we consider midpoint quadrature rules to compute the averaged values. Using $u_{i,j}^n \approx \tilde u_{i,j}^n \approx u(x_i,y_j,t^n)$, $F_{i+1/2,j} \approx f_{i+1/2,j} \approx \\ f(u(x_{i+1/2},y_j,t^{n+1/2}))$, $G_{i,j+1/2} \approx g_{i,j+1/2} \approx g(u(x_i,y_{j+1/2},t^{n+1/2}))$, the numerical scheme then takes the form
\begin{align}
    \label{numericalschemeGeneral}
        u_{i,j}^{n+1} - u_{i,j}^n + \frac{\tau}{h} (F_{i+1/2,j} - F_{i-1/2,j}) + \frac{\tau}{h} (G_{i,j+1/2} - G_{i,j-1/2})  = 0 \,.
\end{align}

To construct the numerical fluxes, we follow standard approaches \cite{shu_essentially_1998,leveque_finite_2004,duraisamy_implicit_2007} of proposing approximations
\rv{
\begin{equation}
    \label{face-values}\nonumber
    u_{i+1/2,j}^{n+1/2,\pm} \approx \lim \limits_{x\rightarrow x_{i+1/2}^{\pm}} u(x,y_j,t^{n+1/2}) \,, \quad
    u_{i,j+1/2}^{n+1/2,\pm} \approx \lim \limits_{y\rightarrow y_{j+1/2}^{\pm}} u(x_i,y,t^{n+1/2}) \,,
\end{equation}
that are} obtained by linear reconstructions in related finite volumes, namely $u_{i+1/2,j}^{n+1/2,-}$ in $I_{i,j}$ and $u_{i+1/2,j}^{n+1/2,+}$ in $I_{i+1,j}$, and analogously for $u_{i,j+1/2}^{n+1/2,\pm}$. 
Having such values and a chosen numerical flux function $H=H(u^-,u^+)$ \cite{shu_essentially_1998,leveque_finite_2004}, one defines
\begin{align}
    \label{GodunovF}
    F_{i+1/2,j} = 
H(u_{i+1/2,j}^{n+1/2,-}, u_{i+1/2,j}^{n+1/2,+})
    \,,
\end{align}
and similarly for the flux $G_{i,j+1/2}$, 
\begin{align}
    \label{GodunovG}
    G_{i,j+1/2} = H( u_{i,j+1/2}^{n+1/2,-}, u_{i,j+1/2}^{n+1/2,+}) \,.
\end{align}
For our purposes, we choose the Godunov flux $H$ defined by
\begin{align}
    \label{GodunovH}
    H( u^{-}, u^{+}) = \begin{cases}
        \min \limits_{u^{-} \leq u \leq u^{+}} h(u) \hspace{1cm} \text{if } u^{-} \leq u^{+} \\
        \max \limits_{u^{+} \leq u \leq u^{-}} h(u) \hspace{1cm} \text{if } u^{-} > u^{+}
    \end{cases} \,,
\end{align}
where $h=f$ for \eqref{GodunovF} and $h=g$ for \eqref{GodunovG}.

Defining in the first step the notation of finite volume methods for the nonlinear scalar equation, in the second step we give additional details for the two dimensional linear advection equation with variable velocity. Here, the flux function $\boldsymbol{f}$ in (\ref{2D_general1}) is formally represented by $\boldsymbol{f}=\Vec{v}u$ with a given space dependent  velocity field $\Vec{v}= \Vec{v}(x,y)=(v(x,y),w(x,y))$. The flux function $\boldsymbol{f}$ now depends on the space variables that we consider only for the linear advection, so the problem \eqref{2D_general1} takes the form
\begin{equation}
    \label{Linear2D}
    \partial_{t} u +  \partial_x (vu) +  \partial_y (wu) = 0 \,.
\end{equation}

In this case, the numerical fluxes are defined as follows,
\begin{align}
    \label{fluxAdvection}
    F_{i+1/2,j} &= v^+_{i+1/2,j}{u}^{n+1/2,-}_{i+1/2,j} + v^-_{i+1/2,j}{u}^{n+1/2,+}_{i+1/2,j} \,, \nonumber \\
    G_{i,j+1/2} &= w^+_{i,j+1/2}{u}^{n+1/2,-}_{i,j+1/2} + w^-_{i,j+1/2}{u}^{n+1/2,+}_{i,j+1/2} \,,
\end{align}
where
\begin{equation}
    \label{AverageVel}
        v_{i\pm 1/2,j} := v(x_{i\pm 1/2},y_j) \,, \quad 
        w_{i,j\pm 1/2} :=  w(x_i,y_{j\pm 1/2}) \,,
\end{equation} 
and we applied the splitting of the velocity vector field,
\begin{align}
    \label{velocitysplit}
    v=v^+ + v^-, \hspace{1cm} & \text{with  } \, v^+ \coloneqq \max(0,v), \,\, v^{-} \coloneqq \min (0,v) \,, \nonumber \\
    w=w^+ + w^-, \hspace{1cm} & \text{with  } \, w^+ \coloneqq \max(0,w), \,\, w^{-} \coloneqq \min (0,w) \,.
\end{align}
In summary, the numerical scheme \eqref{numericalschemeGeneral} will transform into the form 
\begin{align}
    \label{numericalscheme}
         u_{i,j}^{n+1} - {u}_{i,j}^n  &+\frac{\tau}{h} (v^+_{i+1/2,j}{u}^{n+1/2,-}_{i+1/2,j}-v^+_{i-1/2,j}{u}^{n+1/2,-}_{i-1/2,j}) +\frac{\tau}{h} (w^+_{i,j+1/2}{u}^{n+1/2,-}_{i,j+1/2}-w^+_{i,j-1/2}{u}^{n+1/2,-}_{i,j-1/2}) \nonumber\\
         & + \frac{\tau}{h} (v^-_{i+1/2,j}{u}^{n+1/2,+}_{i+1/2,j} - v^-_{i-1/2,j}{u}^{n+1/2,+}_{i-1/2,j} ) + \frac{\tau}{h} (w^-_{i,j+1/2}{u}^{n+1/2,+}_{i,j+1/2} - w^-_{i,j-1/2}{u}^{n+1/2,+}_{i,j-1/2}) = 0 \,.
\end{align}

\rv{Before we define the face values $u_{i\pm1/2,j}^{n+1/2}$ and $u_{i,j\pm1/2}^{n+1/2}$ to be used in the second order accurate compact implicit scheme with \eqref{numericalschemeGeneral} or \eqref{numericalscheme}, we discuss the fast  iterative solvers in the case of the first order accurate approximation that can be used to solve the resulting algebraic systems of equations.
}

\section{Fast sweeping methods}
\label{sec-fsm}

\rv{
To introduce the subject of a special class of iterative solvers to solve the algebraic systems of equations that are obtained by implicit time discretization of the considered hyperbolic equations, we begin with}
the simplest first order accurate approximations 
\begin{align}
\label{firstorder}
    u_{i+1/2,j}^{n+1/2,-}=u_{i,j}^{n+1}\,, \quad u_{i+1/2,j}^{n+1/2,+}=u_{i+1,j}^{n+1} \,, \quad
    u_{i,j+1/2}^{n+1/2,-}=u_{i,j}^{n+1}\,, \quad u_{i,j+1/2}^{n+1/2,+}=u_{i,j+1}^{n+1} \,.
\end{align}

The discretization scheme \eqref{numericalscheme} with \eqref{firstorder} creates a system of algebraic equations that can be solved iteratively using the so called fast sweeping method \cite{tsaiFastSweepingAlgorithms2003,zhao2005fast,chenLaxFriedrichsFast2013}, where each ``sweeping iteration'' is given by standard Gauss-Seidel iteration applied sequentially with alternating index directions across the computational domain. In particular, we use four different $"sweeps"$ from each corner of a rectangular domain: 
\begin{align}
\label{sweeps}
    i = 1, ..., M \ &, \ \ j = 1, ..., M \nonumber \\
    i = M, ..., 1 \ &, \ \ j = 1, ..., M \nonumber \\
    i = M, ..., 1 \ &, \ \ j = M, ..., 1 \nonumber \\ 
    i = 1, ..., M \ &, \ \ j = M, ..., 1  
\end{align}
where each $sweep$ represents the Gauss-Seidel method (GS) with the given ordering as explained next. 
\rv{
Note that the notation in \eqref{sweeps} shall be understood as the outer loop for the index $i$ and the inner loop for the index $j$.
}

\rv{
The Gauss-Seidel method in the case of the first order accurate scheme consists of iterative procedure to obtain $u_{i,j}^{n+1,k+1}\approx u_{i,j}^{n+1} $ by expressing it explicitly from the equation for each $i$ and $j$ written in the form
\begin{align}
    \label{numericalschemeGS}\nonumber
         u_{i,j}^{n+1,k+1}\left(1 +\frac{\tau}{h} \left(v^+_{i+1/2,j} + w^+_{i,j+1/2} - v^-_{i-1/2,j} - w^-_{i,j-1/2}\right)\right) = \\ {u}_{i,j}^n + \frac{\tau}{h} \left(v^+_{i-1/2,j}{u}^{n+1,*}_{i-1,j} + w^+_{i,j-1/2}{u}^{n+1,*}_{i,j-1} - v^-_{i+1/2,j}{u}^{n+1,*}_{i+1,j}  -w^-_{i,j+1/2}{u}^{n+1,*}_{i,j+1}\right) \,,
\end{align}
where $*=k+1$ if the value with this index is available or $*=k$ otherwise. Typically, $u_{i,j}^{n+1,0}=u_{i,j}^n$.  In simple cases like velocities with components not changing their signs, one of the four different sweeps in \eqref{sweeps} is enough to obtain $u_{i,j}^{n+1,1}=u_{i,j}^{n+1}$. For instance, if $v>0$ and $w>0$ then for $i=1,2,\ldots$ and $j=1,2\ldots$ one has always $*=k+1$ for all nonzero terms in the right hand side of \eqref{numericalschemeGS}. Clearly, the values $u_{0,j}^{n+1}$ and $u_{i,0}^{n+1}$ must be defined in this case by boundary conditions. Similar considerations are valid for the nonlinear case \eqref{numericalschemeGeneral} when the sweeps must be realized with nonlinear Gauss-Seidel method \cite{ortega1970iterative} applying, e.g., Newton solver to obtain each value $u_{i,j}^{n+1,k+1}$.
}

\rv{
To understand the convergence of the fast sweeping method for a general velocity field or even for the nonlinear case \eqref{numericalschemeGeneral}, one can refer to a large amount of its applications and variants in the literature. The popularity of the method began with its applications to time independent Hamilton-Jacobi equations \cite{tsaiFastSweepingAlgorithms2003} and especially to eikonal equation \cite{zhao2005fast} . These nonlinear PDEs share the property that the solution can be determined along characteristic curves that is fully respected by the upwind type of space discretization in the fast sweeping method. A complex shape of such curves can increase the number of necessary sweeps, but it is expected that a finite number of iterations (that do not change significantly with a grid refinement) is sufficient to obtain a convergent solution to the algebraic system \cite{tsaiFastSweepingAlgorithms2003,zhao2005fast}, at least at the order of the error introduced by numerical discretization.
}

\rv{
In the context of scalar hyperbolic equations, as we demonstrate later for the chosen relatively simple examples, we can confirm such behavior for numerical solutions obtained with the first order method or our high-resolution compact implicit method. Namely, the expected order of accuracy with respect to refinements of discretization steps is obtained already when the algebraic system is solved iteratively with four sweeps \eqref{sweeps} using standard Gauss-Seidel iterations. Moreover, additional 4 or 8 sweeps decrease the norms of residuals and improve further the accuracy of numerical solutions. 
}

\rv{
The fast sweeping method with standard nonlinear Gauss-Seidel iterations is used for steady state conservation laws employing high order Lax–Friedrichs WENO numerical fluxes in \cite{chenLaxFriedrichsFast2013,chenLaxFriedrichsMultigrid2015} or even with fixed-point sweeping iterations in \cite{wuHighOrderFixedPoint2016} that are based on linear Gauss-Seidel iterations. In those works, the number of GS iterations can grow significantly with the grid refinement. 
}

\rv{
On the other hand, several techniques are based on the fast sweeping strategy that try to preserve its efficiency for a fixed number of iterations like the fast sweeping method with a selection principle based on Rankine-Hugoniot conditions  \cite{engquistFastSweepingMethods2013}, the fractional time step approach \cite{lozano_implicit_2021} or Lower-Upper Symmetric Gauss-Seidel method \cite{yoonLowerupperSymmetricGaussSeidelMethod1988} that is used also as a preconditioner for more involved algebraic solvers \cite{luoFastMatrixfreeImplicit1998a,bocharovImplicitMethodSolution2020}. All these methods are first order accurate in time as they are used typically to obtain stationary solutions, but they suggest that the fast sweeping technique can be well suited even for more complex problems of conservation laws. Therefore, we think that the compact implicit high-resolution method coupled with the fast sweeping method as given in next sections is a promising method to obtain non-oscillatory numerical solutions for time dependent hyperbolic problems with significantly improved accuracy with respect to the first order accurate method.
}

\section{The second order accurate parametric compact implicit scheme}
\label{sec-compact}

In this section, we present the parametric form of the second order accurate scheme in time and space that is later used also in the framework of high-resolution schemes. As discussed above, one has to define the \vr{one sided} approximations $u_{i+1/2,j}^{n+1/2,\pm}$ and $u_{i,j+1/2}^{n+1/2,\pm}$ \vr{to be used in numerical fluxes \eqref{GodunovF} and \eqref{GodunovG}}. For simplicity of the derivation, we will treat in detail only the derivation for the values $u_{i+1/2,j}^{n+1/2,-}$ and $u_{i,j+1/2}^{n+1/2,-}$. Similarly, one can handle the values $u_{i+1/2,j}^{n+1/2,+}$ and $u_{i,j+1/2}^{n+1/2,+}$ for which we will present only the final definitions.

To approximate the values $u_{i+1/2,j}^{n+1/2,-}$ and $u_{i,j+1/2}^{n+1/2,-}$ up to the second order of accuracy \cite{duraisamy_implicit_2007}, we apply a finite Taylor series,
\begin{align}
    \label{taylorX}
    u(x_{i+1/2},y_{j},t^{n+1/2}) = u(x_{i},y_{j},t^{n+1}) + \frac{h}{2} \partial_x u(x_{i}, y_{j}, t^{n+1})  - \frac{\tau}{2}  \partial_t u(x_{i}, y_{j}, t^{n+1}) 
     + \mathcal{O}(h^2, \tau^{2}) \,,
\end{align}   
and
\begin{align}
    \label{taylorY}     
    u(x_{i},y_{j+1/2},t^{n+1/2}) = u(x_{i},y_{j},t^{n+1}) + \frac{h}{2} \partial_y u(x_{i}, y_{j}, t^{n+1})  - \frac{\tau}{2} \partial_t u(x_{i}, y_{j}, t^{n+1}) 
     + \mathcal{O}(h^2, \tau^{2}) \,.
\end{align}
Similarly to \cite{duraisamy_implicit_2007}, we do not replace the time derivative using the Lax-Wendroff procedure, but we approximate \eqref{taylorX} and \eqref{taylorY} directly. Formally, the desired numerical values would then be defined as
\begin{align}
        u_{i+1/2,j}^{n+1/2,-} = u_{i,j}^{n+1} + \frac{h}{2} \partial_x u_{i,j}^{n+1} - \frac{\tau}{2} \partial_t u_{i,j}^{n+1}  
        \,, \quad 
        \label{desired_numerical_values}
        u_{i,j+1/2}^{n+1/2,-} = u_{i,j}^{n+1} + \frac{h}{2} \partial_y u_{i,j}^{n+1} - \frac{\tau}{2} \partial_t u_{i,j}^{n+1} \, ,
\end{align}
where the approximations $\partial_x u_{i,j}^{n+1}$,\,$\partial_y u_{i,j}^{n+1}$ and $\partial_t u_{i,j}^{n+1}$ are expressed differently to \cite{duraisamy_implicit_2007}, namely as a linear combination of two different approximation choices \cite{frolkovic2023high} using parameters $\omega \in [0,1]$,
\begin{align}
    \label{x_approx}
        u_{i+1/2,j}^{n+1/2,-} =  u_{i,j}^{n+1} & + \frac{1}{2}\big(\omega^{x,-} (u_{i,j}^{n+1} - u_{i-1,j}^{n+1}) + (1-\omega^{x,-}) (u_{i+1,j}^{n+1} - u_{i,j}^{n+1})\big) \nonumber \\
        &- \frac{1}{2}\big(\omega^{x,-} (u_{i,j}^{n+1} - u_{i,j}^{n}) + (1-\omega^{x,-}) (u_{i+1,j}^{n+1} - u_{i+1,j}^{n}) \big) \nonumber \\
        = u_{i,j}^{n+1} &- \frac{1}{2}\big(\omega^{x,-}(u_{i-1,j}^{n+1} - u_{i,j}^{n}) + (1-\omega^{x,-})(u_{i,j}^{n+1} - u_{i+1,j}^{n})\big) \,.
\end{align}

Notice that we have managed to obtain a compact stencil as the values $u_{i+1,j}^{n+1}$ in \eqref{x_approx} are canceled. This is also the substantial difference to the approach given in \cite{duraisamy_implicit_2007}. Two particular choices of parameters $\omega^{x,-} \in [0,1]$  do not use the full stencil in \eqref{x_approx}  that we use later conveniently for ENO approximation \cite{shu_essentially_1998}. In particular, with the value $\omega^{x,-}=0$ we create in \eqref{x_approx}  a \emph{"central"} kind of discretization and for $\omega^{x,-}=1$ the \emph{"upwind"} one. 

Very similarly, one obtains the numerical approximations of $u_{i,j+1/2}^{n+1/2,-}$, $u_{i+1/2,j}^{n+1/2,+}$ and $u_{i,j+1/2}^{n+1/2,+}$ as
\begin{align}
        u_{i,j+1/2}^{n+1/2,-} &= u_{i,j}^{n+1} - \frac{1}{2}\big(\omega^{y,-} (u_{i,j-1}^{n+1} - u_{i,j}^{n}) + (1-\omega^{y,-})(u_{i,j}^{n+1} - u_{i,j+1}^{n})\big) \,,   \nonumber \\  \label{u_approx_negat}
        u_{i+1/2,j}^{n+1/2,+} &= u_{i+1,j}^{n+1} - \frac{1}{2}\big(\omega^{x,+} (u_{i+2,j}^{n+1} - u_{i+1,j}^{n}) + (1-\omega^{x,+})(u_{i+1,j}^{n+1} - u_{i,j}^{n})\big) \,, \\ \nonumber
        u_{i,j+1/2}^{n+1/2,+} &= u_{i,j+1}^{n+1} - \frac{1}{2}\big(\omega^{y,+} (u_{i,j+2}^{n+1} - u_{i,j+1}^{n}) + (1-\omega^{y,+})(u_{i,j+1}^{n+1} - u_{i,j}^{n})\big) \,.
\end{align}
For convenience, we gather all parameters by $\boldsymbol{\omega}$, i.e., ${\omega}^{x,\pm} \in \boldsymbol{\omega}$ and ${\omega}^{y,\pm} \in  \boldsymbol{\omega}$. 

We note that the approach presented here is used in \cite{zakovaHigherOrderCompact2024a} to extend the accuracy further in the 1D case up to the third order for smooth solutions with the help of the fractional step method. In that case, one has to use the Lax-Wendroff procedure with mixed derivatives of the solution \cite{zorio_approximate_2017,frolkovic2018semi,carrillo2021lax} to preserve the compactness of the scheme.

To finalize the second order accurate compact implicit scheme, we express the numerical fluxes in \eqref{numericalschemeGeneral} using \eqref{GodunovF} and \eqref{GodunovG} with the numerical approximations $u_{i+1/2,j}^{n+1/2,\pm}, u_{i,j+1/2}^{n+1/2,\pm}$ from \eqref{x_approx} - \eqref{u_approx_negat}. 

\section{High-resolution schemes \rv{with time limiter}}
\label{sec-hr}

For discontinuous initial conditions or when shocks are present in the solution of (\ref{2D_general}), unphysical oscillations may occur in numerical solutions if the numerical methods from the previous section with fixed parameters are used.
To avoid such oscillations, we choose the values of the parameters $\boldsymbol{\omega}$ in each numerical flux differently, i.e., depending on the numerical solution \cite{shu_essentially_1998, frolkovic2023high}.

Let $\omega^{x,\pm}=\omega^{x,\pm}_{i,j} \in [0,1]$ and $\omega^{y,\pm}=\omega^{y,\pm}_{i,j}  \in [0,1]$ in $\boldsymbol{\omega}$ for each $I_{i,j}$ in \eqref{x_approx} - \eqref{u_approx_negat} be the free parameters that we want to determine. 
As is known from the literature \cite{duraisamy_implicit_2007, arbogast2020third, puppo_quinpi_2022, frolkovic2023high,puppo2023quinpi}, unphysical oscillations can occur not only due to an inappropriately fixed choice of space reconstruction, but also due to the fixed time reconstruction. Therefore, similarly to \cite{frolkovic2023high}, we add another numerical parameters gathered in $\boldsymbol{l}=({l}^{x,\pm}_{i,j}, {l}^{y,\pm}_{i,j}) $, \vr{with ${l}^{x,\pm}_{i,j}, {l}^{y,\pm}_{i,j} \in [0,1]$,} which, if necessary, can limit the second order space reconstruction towards the first order form of the scheme to produce numerical solutions free of unphysical oscillations. 
Notice that the values $\boldsymbol{\omega}$ and $\boldsymbol{l}$ also change at each time step, which we do not emphasize in the notation. 

In summary, the approximation \eqref{x_approx} - \eqref{u_approx_negat} will transform into the form
\begin{align}
    \label{2D_ENO_pozit}
         u_{i+1/2,j}^{n+1/2,-} &= u_{i,j}^{n+1} - \frac{l^{x,-}_{i,j}}{2}\big(\omega^{x,-}_{i,j} (u_{i-1,j}^{n+1} - 
            u_{i,j}^{n}) + (1-\omega^{x,-}_{i,j})(u_{i,j}^{n+1} - u_{i+1,j}^{n})\big) \,, \nonumber \\
         u_{i,j+1/2}^{n+1/2,-} &= u_{i,j}^{n+1} - \frac{l^{y,-}_{i,j}}{2}\big(\omega^{y,-}_{i,j} (u_{i,j-1}^{n+1} - 
            u_{i,j}^{n}) + (1-\omega^{y,-}_{i,j})(u_{i,j}^{n+1} - u_{i,j+1}^{n})\big)\, ,
\end{align}
and for (\ref{u_approx_negat})
\begin{align}
    \label{2D_ENO_negat}
        u_{i+1/2,j}^{n+1/2,+} &= u_{i+1,j}^{n+1} - \frac{l^{x,+}_{i+1,j}}{2}\big(\omega^{x,+}_{i+1,j} (u_{i+2,j}^{n+1} - u_{i+1,j}^{n}) + (1-\omega^{x,+}_{i+1,j})(u_{i+1,j}^{n+1} - u_{i,j}^{n})\big) \,, \nonumber\\
        u_{i,j+1/2}^{n+1/2,+} &= u_{i,j+1}^{n+1} - \frac{l^{y,+}_{i,j+1}}{2}\big(\omega^{y,+}_{i,j+1} (u_{i,j+2}^{n+1} - u_{i,j+1}^{n}) + (1-\omega^{y,+}_{i,j+1})(u_{i,j+1}^{n+1} - u_{i,j}^{n})\big) \,.
\end{align}
The high-resolution scheme is then defined by \eqref{numericalschemeGeneral} - \eqref{GodunovH} with the face values given in \eqref{2D_ENO_pozit} - \eqref{2D_ENO_negat} \rv{
for which the parameters in $\boldsymbol{\omega}$ and $\boldsymbol{l}$ will depend on the numerical solution to suppress its oscillatory behavior. The parameters in $\boldsymbol{\omega}$ will be used to choose appropriate space reconstruction in \eqref{x_approx} - \eqref{u_approx_negat} and the parameters in $\boldsymbol{l}$ to limit, eventually, the fixed higher order time discretization towards the first order accuracy.
}

Concerning the space reconstruction, we present two standard choices - the simplest Essentially Non-Oscillatory (ENO) method and a variant of the Weighted ENO (WENO) method. The value of $\boldsymbol{\omega}$ for the ENO method will depend on the following ratios in $\boldsymbol{r}=({r}^{x,\pm}_{i,j}, {r}^{y,\pm}_{i,j})$,
\begin{align}
    \label{r}
    r^{x,-}_{i,j}=\frac{ u_{i-1,j}^{n+1} - u_{i,j}^{n} }{u_{i,j}^{n+1} - u_{i+1,j}^{n}} \,, &\hspace{1.5cm} r^{x,+}_{i,j}=\frac{ u_{i+1,j}^{n+1} - u_{i,j}^{n} }{u_{i,j}^{n+1} - u_{i-1,j}^{n}} \,,\nonumber \\
    r^{y,-}_{i,j}=\frac{ u_{i,j-1}^{n+1} - u_{i,j}^{n} }{u_{i,j}^{n+1} - u_{i,j+1}^{n}} \,, &\hspace{1.5cm} r^{y,+}_{i,j}=\frac{ u_{i,j+1}^{n+1} - u_{i,j}^{n} }{u_{i,j}^{n+1} - u_{i,j-1}^{n}} \,,
\end{align}
while the WENO method will depend on the nominators and denominators in the definitions of $\boldsymbol{r}$.

In particular, we can define the Essentially Non-Oscillatory (ENO) scheme \cite{shu_essentially_1998} based on the parameters $\boldsymbol{r}$ simply by
\begin{align}
\label{omegachoice}
    {\omega}^{x,\pm}_{i,j} = 
        \begin{cases}
        1 \hspace{0.5cm} \text{if } |{r}^{x,\pm}_{i,j}| \leq 1 \\
        0 \hspace{0.5cm} \text{otherwise}
         \end{cases} \,, \quad
    {\omega}^{y,\pm}_{i,j} = 
    \begin{cases}
    1 \hspace{0.5cm} \text{if } |{r}^{y,\pm}_{i,j}| \leq 1 \\
    0 \hspace{0.5cm} \text{otherwise}
     \end{cases} \,.
\end{align}

To create a weighted essentially non-oscillatory scheme,
we propose, again, a solution-dependent choice of ${\omega}^{x,\pm}_{i,j}, {\omega}^{y,\pm}_{i,j}$ as found in \cite{shu_essentially_1998} and adapted to our compact implicit scheme. The WENO scheme is used to achieve a smoother choice, compared to the ENO scheme, between the \textit{,,upwind''} and \textit{,,central''} differences in (\ref{x_approx}) - (\ref{u_approx_negat}), so $\omega^{x,\pm}_{i,j}, \omega^{y,\pm}_{i,j}\in[0,1]$.

Let $\bar{\omega}$ be some preferable constant value of $\omega^{x,\pm}_{i,j}, \omega^{y,\pm}_{i,j}$. \rv{In all our numerical experiments we choose $\bar \omega=1/2$.}
Then, using the nominators and denominators in $\boldsymbol{r}$ given in (\ref{r}), we set first the weights
\begin{align}
    \label{WENOweightsPOS}
    a^{x,\pm}_{u} =   \frac{\bar{\omega}}{\left(\epsilon + (u_{i\pm1,j}^{n+1}-u_{i,j}^{n})^2\right)^2} \,,  &\hspace{1cm}
    a^{y,\pm}_{u} =  \frac{\bar{\omega}}{\left(\epsilon + (u_{i,j\pm1}^{n+1}-u_{i,j}^{n})^2\right)^2} \,, \nonumber \\
    a^{x,\pm}_{c} = \frac{1-\bar{\omega}}{\left(\epsilon + (u_{i\mp,j}^{n+1}-u_{i,j}^{n})^2\right)^2} \,, 
    &\hspace{1cm}
    a^{y,\pm}_{c} = \frac{1-\bar{\omega}}{\left(\epsilon + (u_{i,j\mp}^{n+1}-u_{i,j}^{n})^2\right)^2} \,,
\end{align}
with $\epsilon$ being a small number to avoid division by zero. Afterwards, we define the parameters $\omega^{x,\pm}_{i,j}, \omega^{y,\pm}_{i,j}$ by
\begin{equation}
    \label{WENOomega}
    \omega^{x,-}_{i,j} = \frac{a^{x,-}_{u}}{a^{x,-}_{u}+a^{x,-}_{c}} \,, \hspace{0.5cm} \omega^{y,-}_{i,j} = \frac{a^{y,-}_{u}}{a^{y,-}_{u}+a^{y,-}_{c}} \,, \hspace{0.5cm}  \omega^{x,+}_{i,j} = \frac{a^{x,+}_{u}}{a^{x,+}_{u}+a^{x,+}_{c}} \,, \hspace{0.5cm} \omega^{y,+}_{i,j} = \frac{a^{y,+}_{u}}{a^{y,+}_{u}+a^{y,+}_{c}} \,.
\end{equation}

\rv{
Note that the ratios $\boldsymbol{r}$ enable us to express the second order updates of the face values defined in \eqref{2D_ENO_pozit}  and \eqref{2D_ENO_negat} in a multiplicative way, 
\begin{align}
    \label{2D_ENO_pozit-short}
         u_{i+1/2,j}^{n+1/2,-} = u_{i,j}^{n+1} - \frac{l^{x,-}_{i,j}}{2}\frac{\Psi_{i,j}^{x,-}}{r_{i,j}^{x,-}} (u_{i-1,j}^{n+1} - 
            u_{i,j}^{n})  \,, \quad
         u_{i,j+1/2}^{n+1/2,-} = u_{i,j}^{n+1} - \frac{l^{y,-}_{i,j}}{2} \frac{\Psi_{i,j}^{y,-}}{r_{i,j}^{y,-}} (u_{i,j-1}^{n+1} - 
            u_{i,j}^{n}) \,,\\
    \label{2D_ENO_negat-short}
        u_{i+1/2,j}^{n+1/2,+} = u_{i+1,j}^{n+1} - \frac{l^{x,+}_{i+1,j}}{2}\Psi_{i+1,j}^{x,+}(u_{i+1,j}^{n+1} - u_{i,j}^{n}) \,, \quad
        u_{i,j+1/2}^{n+1/2,+} = u_{i,j+1}^{n+1} - \frac{l^{y,+}_{i,j+1}}{2}\Psi_{i,j+1}^{y,+}(u_{i,j+1}^{n+1} - u_{i,j}^{n}) \,,
\end{align}
using
\begin{equation}
    \label{Psi}
    \Psi_{i,j}^{*,\pm} := \omega^{*,\pm}_{i,j} r_{i,j}^{*,\pm} + 1-\omega^{*,\pm}_{i,j} \,,
\end{equation}
where $*=x$ and $*=y$. The values $\Psi$ in \eqref{Psi} are typically used when defining high-resolution schemes using so-called limiters \cite{sweby1984high,spekreijseMultigridSolutionMonotone1987,leveque_finite_2004,kemm_comparative_2011}. In our case they help us to motivate the definitions of the parameters in $\boldsymbol{l}$ as given in the next section. 
}

\rv{
We disclose in advance that we require that the values $\Psi_{i,j}^{*,\pm}$ fulfill for any values of ${r}_{i,j}^{*,\pm} \in R$ the inequalities \cite{spekreijseMultigridSolutionMonotone1987,kemm_comparative_2011,zhangReviewTVDSchemes2015,woodfieldNewLimiterRegions2024}
\begin{equation}
    \label{psi-inequalities}
    -1 \le \Psi_{i,j}^{*,\pm} \le 2 \,, \quad i,j,=1,2,\ldots,I \,.
\end{equation}
One can show by a simple inspection that these inequalities are fulfilled by the ENO approximation \eqref{omegachoice} using \eqref{Psi}.
}

\rv{\section{Positive coefficient scheme}}
\label{sec-positive}

In this section, we formulate nonlinear conditions that \rv{the parameters  $\boldsymbol l$ of the time limiter} in (\ref{2D_ENO_pozit}) and (\ref{2D_ENO_negat}) must satisfy so that the high-resolution scheme (\ref{numericalscheme}) for the linear advection equation (\ref{Linear2D}) does not produce unphysical oscillations for divergence free velocity. For that purpose, we formulate the conditions when it can be written in the form \rv{of positive coefficient scheme} \cite{spekreijseMultigridSolutionMonotone1987,jamesonPositiveSchemesShock1995,woodfieldNewLimiterRegions2024},
\begin{align}
\label{positive}
    P_{i,j} (u_{i,j}^{n+1}-u_{i,j}^{n}) &+ P_{i+1,j}( u_{i,j}^{n+1}- u_{i+1,j}^{n+1})+P_{i-1,j}( u_{i,j}^{n+1}- u_{i-1,j}^{n+1}) \nonumber \\
    &+ P_{i,j+1}( u_{i,j}^{n+1}- u_{i,j+1}^{n+1})+P_{i,j-1}( u_{i,j}^{n+1}- u_{i,j-1}^{n+1}) = 0 \,,
\end{align}
where all five coefficients $P$ are non-negative. Note that these coefficients are different for each discrete equation. Clearly, if (\ref{positive}) is valid, a local discrete minimum and maximum principle is fulfilled as $u_{i,j}^{n+1}$ is a convex combination of all neighboring values. To express (\ref{numericalscheme}) \rv{with \eqref{2D_ENO_pozit} and \eqref{2D_ENO_negat}} in the form (\ref{positive}), all coefficients $P$ must depend on unknown values of the numerical solution.

The first important assumption is that the discrete values of velocities must fulfill the discrete condition of incompressibility $\nabla \cdot \vec{v} = 0$, i.e., 
\begin{align}
\label{incompr}
    C_{i+1/2,j}^{+}-C_{i-1/2,j}^{-}+C_{i+1/2,j}^{-}-C_{i-1/2,j}^{+}+C_{i,j+1/2}^{+}-C_{i,j-1/2}^{-}+C_{i,j+1/2}^{-}-C_{i,j-1/2}^{+} = 0 \,,
\end{align}
\rv{
with
\begin{align}
    C^{+}_{i\pm1/2,j}=\frac{\tau}{h}v^+_{i\pm1/2,j}\,, \quad C^{-}_{i\pm1/2,j}=\frac{\tau}{h}v^-_{i\pm1/2,j}\,, \quad
    C^{+}_{i,j\pm1/2}=\frac{\tau}{h}w^+_{i,j\pm1/2}\,, \quad C^{-}_{i,j\pm1/2}=\frac{\tau}{h}w^-_{i,j\pm1/2}\,.
\end{align}
Subtracting (\ref{incompr}) multiplied by $u_{i,j}^{n+1}$ from (\ref{numericalscheme}) we obtain its equivalent form,
\begin{align}
    \label{numericalscheme-uij-subtracted}
         u_{i,j}^{n+1} - {u}_{i,j}^n  &+C^{+}_{i+1/2,j}({u}^{n+1/2,-}_{i+1/2,j} -u_{i,j}^{n+1})-C^{+}_{i-1/2,j}({u}^{n+1/2,-}_{i-1/2,j}-u_{i,j}^{n+1}) 
         \nonumber\\
         & + C^{-}_{i+1/2,j}({u}^{n+1/2,+}_{i+1/2,j} - u_{i,j}^{n+1}) - C^{-}_{i-1/2,j} ({u}^{n+1/2,+}_{i-1/2,j} - u_{i,j}^{n+1} ) + \ldots = 0\,,
\end{align}
where we have skipped four analogous terms beginning with $C_{i,j\pm1/2}^{\pm}$.}
One can easily show that the first order scheme can be transferred to the form (\ref{positive}) using \eqref{numericalscheme-uij-subtracted} with (\ref{firstorder}), when 
\begin{align}
    P_{i,j}=1, \quad P_{i+1,j}=-C_{i+1/2,j}^{-},  \quad P_{i-1,j}= C_{i-1/2,j}^{+},  \quad P_{i,j+1}=-C_{i,j+1/2}^{-}, 
     \quad P_{i,j-1}=C_{i,j-1/2}^{+} \,.
\end{align}
\rv{Clearly, all coefficients $P$ in above are nonnegative.} Consequently, if the parameters in \rv{time limiters $\boldsymbol{l}$ are set to zero in \eqref{2D_ENO_pozit} and \eqref{2D_ENO_negat} for the high-resolution scheme} (\ref{numericalscheme}), we obtain the first order form of the scheme that produces numerical solutions free of unphysical oscillations.

\rv{Next, we derive an equivalent form of the high-resolution scheme \eqref{numericalscheme} with \eqref{2D_ENO_pozit} and \eqref{2D_ENO_negat} for which the definitions of coefficients $P$ in \eqref{positive} will follow. First, to make the notations shorter, we introduce,
\begin{eqnarray}
    \label{AD}
    D_{i,j}^{*,\pm} := \frac{1}{2} \frac{\Psi_{i,j}^{*,\pm}}{r_{i,j}^{*,\pm}} \,, \quad
    A_{i\pm1,j} := \frac{1}{2} l_{i\pm1,j}^{x,\pm} \Psi_{i\pm1,j}^{x,\pm} \,, \quad
    A_{i,j\pm1} := \frac{1}{2} l_{i,j\pm1}^{y,\pm} \Psi_{i,j\pm1}^{y,\pm} \,,
\end{eqnarray}
where $*=x$ and $*=y$. Consequently, we obtain from \eqref{2D_ENO_pozit-short} - \eqref{2D_ENO_negat-short} 
}
\begin{align}
    \label{uhr1}
    u_{i+1/2,j}^{n+1/2,-} - u_{i,j}^{n+1} & = D_{i,j}^{x,-} l_{i,j}^{x,-} (
            u_{i,j}^{n}-u_{i-1,j}^{n+1}) \,, \\[1ex]
                \label{uhr2}
     u_{i+1/2,j}^{n+1/2,+} - u_{i,j}^{n+1} & =  - \Big(u_{i,j}^{n+1} - u_{i+1,j}^{n+1} - A_{i+1,j} ( u_{i,j}^{n} - u_{i+1,j}^{n+1}) \Big) \,,
\end{align}
and analogously for other terms with $u_{i-1/2,j}^{n+1/2,\pm}$ and $u_{i,j\pm 1/2}^{n+1/2,\pm}$. \rv{ Using \eqref{uhr1} - \eqref{uhr2}, 
we obtain the following form of \eqref{numericalscheme-uij-subtracted},}
\begin{align}
    \label{hr-cond2}
    u_{i,j}^{n+1}-u^n_{i,j} & + C_{i+1/2,j}^{+} l^{x,-}_{i,j} D^{x,-}_{i,j} ( 
            u_{i,j}^{n} - u_{i-1,j}^{n+1})  
            - \, C_{i+1/2,j}^{-}\left( u_{i,j}^{n+1} - u_{i+1,j}^{n+1}  - A_{i+1,j}   (u_{i,j}^{n} - u_{i+1,j}^{n+1}) \right)  \nonumber \\[1ex] 
        & -\, C_{i-1/2,j}^{-} l^{x,+}_{i,j}  ( 
            u_{i,j}^{n}-u_{i+1,j}^{n+1}) 
            + \, C_{i-1/2,j}^{+}\left( u_{i,j}^{n+1} - u_{i-1,j}^{n+1}  - A_{i-1,j}  (u_{i,j}^{n} - u_{i-1,j}^{n+1}) \right) + \ldots = 0 \,,
\end{align}
where we have skipped analogous four terms that should appear after $C^{\pm}_{i,j\pm1/2}$. Taking into account that
$$
u_{i,j}^n - u_{i\pm1,j}^{n+1} = u_{i,j}^{n+1} - u_{i\pm1,j}^{n+1} - (u_{i,j}^{n+1} - u_{i,j}^{n}) \,,
$$
and \rv{analogously for $u_{i,j}^n - u_{i,j\pm1}^{n+1}$}, the equations (\ref{hr-cond2}) can be written in the form (\ref{positive}) with the coefficients
\begin{align}
\label{Px}
    P_{i-1,j} = C_{i+1/2,j}^{+} &l_{i,j}^{x,-} D_{i,j}^{x,-} + C_{i-1/2,j}^{+} (1 - A_{i-1,j}) \,, \quad 
    P_{i+1,j} = -C_{i-1/2,j}^{-} l_{i,j}^{x,+} D_{i,j}^{x,+} - C_{i+1/2,j}^{-} (1 - A_{i+1,j}) \,, \\[1ex]
\label{Py}
    P_{i,j-1} = C_{i,j+1/2}^{+} &l_{i,j}^{y,-} D_{i,j}^{y,-} + C_{i,j-1/2}^{+} (1 - A_{i,j-1}) \,, \quad 
    P_{i,j+1} = -C_{i,j-1/2}^{-} l_{i,j}^{y,+} D_{i,j}^{y,+} - C_{i,j+1/2}^{-} (1 - A_{i,j+1}) \,, \\[1ex]
\label{Pij}
    P_{i,j} = 1 &-
    C_{i+1/2,j}^{+}  l_{i,j}^{x,-} D_{i,j}^{x,-} + C_{i-1/2,j}^{+} A_{i-1,j} +
    C_{i-1/2,j}^{-}  l_{i,j}^{x,+} D_{i,j}^{x,+} - \left. C_{i+1/2,j}^{-} A_{i+1,j} \right.  \nonumber \\[1ex] 
    &- C_{i,j+1/2}^{+}  l_{i,j}^{y,-} D_{i,j}^{y,-} + C_{i,j-1/2}^{+} A_{i,j-1} +
    C_{i,j-1/2}^{-}  l_{i,j}^{y,+} D_{i,j}^{y,+} - C_{i,j+1/2}^{-} A_{i,j+1} \,. 
\end{align}

\rv{
To have nonnegative coefficients in \eqref{Px} and \eqref{Py}, we have to require that 
\begin{equation}
    \label{Aless1}
    A_{i\pm1,j}\le1\,, \quad A_{i,j\pm1}\le1 \,,
\end{equation}
that is clearly fulfilled for the assumptions in \eqref{psi-inequalities} due to \eqref{AD}.
Consequently, one can always find the time limiter parameters in $\boldsymbol{l}$ such that
\begin{equation}
    \label{l1}
    |C_{i\pm1/2,j}^{\pm}| l_{i,j}^{x,\mp} D_{i,j}^{x,\mp} \ge - |C_{i\mp1/2,j}^{\pm}| (1 - A_{i\mp1,j}) \,, \quad
    |C_{i,j\pm1/2}^{\pm}| l_{i,j}^{y,\mp} D_{i,j}^{y,\mp} \ge - |C_{i,j\mp1/2}^{\pm}| (1 - A_{i,j\mp1}) \,,
\end{equation}
so the coefficients $P_{i\pm1,j}$ and $P_{i,j\pm1}$ in \eqref{Px} - \eqref{Py} become nonnegative. Note that the inequalities \eqref{l1} are fulfilled trivially if $D_{i,j}^{*,\mp}\ge 0$ and the values in $\boldsymbol{l}$ could be different from $1$ only for $D_{i,j}^{*,\mp}< 0$. In any case, the values $l_{i,j}^{*,\mp}=0$ guarantee the validity of \eqref{l1}. We define the values in $\boldsymbol{l}$ at the end of this section.
}

\rv{
Finally, we have to prove the sign property for $P_{i,j}$ in \eqref{Pij}. Note that this step is unique for the compact implicit scheme as analogous coefficient in \eqref{positive} is typically equal one for the fully explicit or implicit scheme \cite{harten_class_1984, spekreijseMultigridSolutionMonotone1987,jamesonPositiveSchemesShock1995}. First, we define the local (grid) Courant number
\begin{equation}
    \label{CN}
    C_{i,j} := C_{i+1/2,j}^{+} - C_{i-1/2,j}^{-} +  C_{i,j+1/2}^{+} - C_{i,j-1/2}^{-} = - C_{i+1/2,j}^{-} + C_{i-1/2,j}^{+} - C_{i,j+1/2}^{-} + C_{i,j-1/2}^{+} \,.
\end{equation}
Note that
\begin{equation}
    \label{PandP}
    P_{i,j}=1+C_{i,j}-P_{i+1,j}-P_{i-1,j}-P_{i,j+1}-P_{i,j-1} \,,
\end{equation}
and
\begin{equation}
    \label{1plusC}
    1+C_{i,j} = \left(\frac{1}{C_{i,j}} + 1 \right) \left(- C_{i+1/2,j}^{-} + C_{i-1/2,j}^{+} - C_{i,j+1/2}^{-} + C_{i,j-1/2}^{+}\right)  \,,
\end{equation}
so one can reach that $P_{i,j}\ge0$ if
\begin{eqnarray}
    \label{condition-on-l-prepare}
    P_{i\pm1,j} \le \frac{|C_{i\pm1/2,j}^{\mp}|}{C_{i,j}}+|C_{i\pm1/2,j}^{\mp}| \,, \quad P_{i,j\pm1} \le \frac{|C_{i,j\pm1/2}^{\mp}|}{C_{i,j}}+|C_{i,j\pm1/2}^{\mp}| \,,
\end{eqnarray}
or, equivalently,
\begin{eqnarray}
    \label{condition-on-l}
    |C_{i\mp1/2,j}^{\mp}| l_{i,j}^{x,\pm} D_{i,j}^{x,\pm} \le \frac{|C_{i\pm1/2,j}^{\mp}|}{C_{i,j}} + |C_{i\pm1/2,j}^{\mp}| A_{i\pm1,j} \,,
    \quad 
    |C_{i,j\mp1/2}^{\mp}| l_{i,j}^{y,\pm} D_{i,j}^{y,\pm} \le \frac{|C_{i,j\pm1/2}^{\mp}|}{C_{i,j}} + |C_{i,j\pm1/2}^{\mp}| A_{i,j\pm1} \,.
\end{eqnarray}
To have a chance to find the time limiter parameters in $\boldsymbol{l}$ to fulfill the inequalities in \eqref{condition-on-l}, we have to require
\begin{equation}
    \label{Amore}
    A_{i\pm1,j} \ge - \frac{1}{C_{i,j}} \,, \quad A_{i,j\pm1} \ge - \frac{1}{C_{i,j}} \,,
\end{equation}
that can be obtained by the restrictions
\begin{equation}
    \label{Psi2}
    {l}_{i,j}^{x,\pm} \Psi_{i,j}^{x,\pm} \ge - \frac{2}{C_{i\mp1,j}} \,,\quad {l}_{i,j}^{y,\pm} \Psi_{i,j}^{y,\pm} \ge - \frac{2}{C_{i,j\mp1}} \,, \quad i,j=1,2,\ldots,I \,.
\end{equation}
The inequalities \eqref{Psi2} are trivially fulfilled if $\Psi_{i,j}^{*,\pm} \ge 0$.
}

\rv{
Now we are ready to define the values of the time limiter parameters in $\boldsymbol{l}$. In summary, to obtain non-negative coefficients in \eqref{Px} - \eqref{Pij}, we have to fulfill the inequalities \eqref{Aless1} - \eqref{l1} and \eqref{condition-on-l} - \eqref{Amore}. To proceed, we define the values
\begin{equation}
    \label{dfn1-l}
    L_{i,j}^{x,\pm} :=\left(\max \{ 1, -\frac{C_{i\mp1,j} \Psi_{i,j}^{x,\pm}}{2}\} \right)^{-1} \,, \quad
    L_{i,j}^{y,\pm} :=\left(\max \{ 1, -\frac{C_{i,j\mp1} \Psi_{i,j}^{y,\pm}}{2}\} \right)^{-1}\,.
\end{equation}
}

\rv{
Next, we distinguish the following cases. Firstly, if $C_{i\mp1/2,j}^{\mp} D_{i,j}^{x,\pm} = 0$ or $C_{i,j\mp1/2}^{\mp} D_{i,j}^{y,\pm} = 0$, then the inequalities \eqref{l1} and \eqref{condition-on-l} are fulfilled automatically and we set
\begin{equation}
    \label{lequalsL}
    l_{i,j}^{*,\pm} = L_{i,j}^{*,\pm} \,.
\end{equation}
If $|C_{i\mp1/2,j}^{\mp}| D_{i,j}^{*,\pm} > 0$, then we define 
\begin{eqnarray}
    \label{dfn2-l}
    {l}^{x,\pm}_{i,j} = \min \Big{\{} L_{i,j}^{x,\pm} , \frac{C_{i\pm1/2,j}^{\mp}}{C_{i\mp1/2,j}^{\mp} D_{i,j}^{x,\pm} }\left(\frac{1}{C_{i,j}} + A_{i\pm1,j} \right)  \Big{\}} \,, \\
    {l}^{y,\pm}_{i,j} = \min \Big{\{} L_{i,j}^{y,\pm} , \frac{C_{i,j\pm1/2}^{\mp}}{C_{i,j\mp1/2}^{\mp} D_{i,j}^{y,\pm} }\left(\frac{1}{C_{i,j}} + A_{i,j\pm1} \right)  \Big{\}} \,,
\end{eqnarray}
and if $|C_{i\mp1/2,j}^{\mp}| D_{i,j}^{*,\pm} < 0$, then we define 
\begin{eqnarray}
    {l}^{x,\pm}_{i,j} = \min \Big{\{} L_{i,j}^{x,\pm} ,  - \frac{C_{i\pm1/2,j}^{\mp}}{C_{i\mp1/2,j}^{\mp} D_{i,j}^{x,\mp}} (1 - A_{i\pm1,j})  \Big{\}} \,, \\
    \label{dfn3-l}
    {l}^{y,\pm}_{i,j} = \min \Big{\{} L_{i,j}^{y,\pm} ,  - \frac{C_{i,j\pm1/2}^{\mp}}{C_{i,j\mp1/2}^{\mp} D_{i,j}^{x,\mp}} (1 - A_{i,j\pm1})  \Big{\}} \,.
\end{eqnarray}
}

\rv{
In summary, the high-resolution scheme \eqref{numericalscheme} with \eqref{2D_ENO_pozit} and \eqref{2D_ENO_negat} for the linear advection equation \eqref{Linear2D} is the positive coefficient scheme \eqref{positive} if the parameters $\boldsymbol{\omega}$ are defined such that the inequalities \eqref{psi-inequalities} are fulfilled and the parameters in $\boldsymbol{l}$ are defined according to \eqref{dfn1-l} - \eqref{dfn3-l}. Note that the definitions \eqref{dfn2-l} - \eqref{dfn3-l} for ${l}_{i,j}^{x,\pm}$ are coupled with ${l}_{i\pm1,j}^{x,\mp}$ (and analogously for ${l}_{i,j}^{y,\pm}$), but this dependence is resolved iteratively due to the fast sweeping method.
}

\begin{remark}
\rv{In our numerical experiments,} we follow the simplest and robust approach without an attempt to optimize the definitions of $\boldsymbol \omega$ and $\boldsymbol l$ to maximize accuracy.  

Firstly, we require ${l}_{i,j}^{*,\pm}=0$ if ${r}_{i,j}^{*,\pm}<0$. Note that the ratios in $\boldsymbol{r}$ take negative values only near the extrema of numerical solutions. \rv{In particular, we replace the definitions \eqref{lequalsL} by simpler one
\begin{equation}
    \label{Lsimple}
    L_{i,j}^{*,\pm} = \left \{\begin{array}{lr}
       1  & {r}_{i,j}^{*,\pm} \ge 0 \\
       0  & {r}_{i,j}^{*,\pm} < 0
    \end{array}
    \right. \,.
\end{equation}
}

\rv{
Secondly, in this study we consider examples when $C_{i+1/2,j}=C_{i-1/2,j}$ and $C_{i,j+1/2}=C_{i,j-1/2}$, so we need not to consider the dependence of parameters in $\boldsymbol{l}$ on the face Courant numbers $C_{i\pm1,j}^{\pm}$ and  $C_{i,j\pm1}^{\pm}$. Consequently,}  the parameters ${l}^{x,\pm}_{i,j}, {l}^{y,\pm}_{i,j}$ are determined by 
\begin{align}
    \label{l} 
    {l}^{x,\pm}_{i,j} &= \min \Big{\{} 1, \max \big{\{} 0, \left( {\omega}^{x,\pm}_{i,j} + \frac{1-{\omega}^{x,\pm}_{i,j} }{{r}^{x,\pm}_{i,j}}\right)^{-1}\left(\frac{2}{C_{i,j}}+{l}^{x,\pm}_{i\pm1,j}({\omega}^{x,\pm}_{i\pm1,j} {r}^{x,\pm}_{i\pm1,j} + 1-{\omega}^{x,\pm}_{i\pm1,j})\right)\big{\}} \Big{\}} \,, \nonumber \\
    {l}^{y,\pm}_{i,j} &= \min \Big{\{} 1, \max \big{\{} 0, \left({\omega}^{y,\pm}_{i,j} + \frac{1-{\omega}^{y,\pm}_{i,j} }{{r}^{y,\pm}_{i,j}}\right)^{-1}\left(\frac{2}{C_{i,j}}+{l}^{y,\pm}_{i,j\pm1}({\omega}^{y,\pm}_{i,j\pm1} {r}^{y,\pm}_{i,j\pm1} + 1-{\omega}^{y,\pm}_{i,j\pm1})\right)\big{\}} \Big{\}} \,. 
\end{align} 

For the nonlinear scalar equation, \rv{to avoid another nonlinearity due to  the dependence of characteristic speeds on the solution, the local Courant number $C_{i,j}$ in \eqref{l} is replaced by the fixed maximal Courant number $C$  given by}
\begin{equation}\label{C}
    C = \frac{\tau}{h}\max_{u} |f'(u)| +  \frac{\tau}{h} \max_u |g'(u)| \,.
\end{equation}
\end{remark}

\begin{remark}
To evaluate the nonlinear dependence \rv{of parameters in $\boldsymbol{\omega}$ and $\boldsymbol{l}$ on $u_{i,j}^{n+1}$, the following procedure will be used.}
\begin{enumerate}
  \item Calculate an initial estimate of $u_{i,j}^{n+1}$ using \rv{the first order accurate numerical scheme \eqref{numericalschemeGeneral} - \eqref{GodunovH} with \eqref{firstorder} using four GS iterations \eqref{sweeps}. }
  \item \label{Step:CorrectorENO}
  \rv{
  Apply GS iterations \eqref{sweeps} up to desired reduction of the residuals 
  where to compute each value of $u_{i,j}^{n+1,k+1}\approx u_{i,j}^{n+1}$ do the following steps:}
      \begin{enumerate}
      \item Calculate ${r}^{x,\pm}_{i,j}$ and ${r}^{y,\pm}_{i,j}$ from \eqref{r} using the \rv{last available} values of $u_{i,j}^{n+1,k}$. 
      \item Set ${\omega}^{x,\pm}_{i,j}$ and ${\omega}^{y,\pm}_{i,j}$ with (\ref{omegachoice}) based on the calculated $\boldsymbol{r}$ in the step 2.1; afterwards, calculate the values ${l}^{x,\pm}_{i,j}, {l}^{y,\pm}_{i,j}$ from (\ref{l}).
      \item \rv{Update the value $u_{i,j}^{n+1}$ using the high-resolution numerical scheme \eqref{numericalschemeGeneral} - \eqref{GodunovH} and \eqref{2D_ENO_pozit} - \eqref{2D_ENO_negat}} with the obtained values ${\omega}^{x,\pm}_{i,j}$, ${\omega}^{y,\pm}_{i,j}$ and ${l}^{x,\pm}_{i,j}$ and ${l}^{y,\pm}_{i,j}$ in the step 2.2.
    \end{enumerate}
\end{enumerate}

The iterative procedure for using the WENO scheme is analogous to the procedure explained for the ENO scheme. We compute several numerical experiments (see Section \ref{sec-numerical}) for which we investigate, among others, the behavior of the fast sweeping method depending on the number of Gauss-Seidel iterations.
    
\end{remark}

\section{Numerical experiments}
\label{sec-numerical}

We present numerical results of the proposed second order and high-resolution compact implicit finite-volume scheme with the purpose of illustrating their accuracy and stability properties. \vr{The second order unlimited scheme is used as described in Section \ref{sec-compact} to show the expected experimental order of convergence for different fixed values of parameter ${\omega}$. The high-resolution scheme is using the time limiter defined by \eqref{l} and the choice of $\boldsymbol{\omega}$ in \eqref{omegachoice} for the ENO discretization and the choice in \eqref{WENOomega} for the WENO discretization.} We describe the results \vr{for these choices} in two subsections, the first one dealing with a linear advection equation, and the second one dealing with the Burgers equation as a simple representative nonlinear problem. For each of the cases, we present an example with a smooth initial condition, as well as an example with a discontinuous initial condition, employing both the ENO and WENO approximations. Specifically for the Burgers equation, we demonstrate two phenomena, shock and rarefaction waves, separately.

If the exact solution is known for the chosen example, it is used for boundary conditions, and the discrete $L_1$ norm ($E$) of the error is calculated as follows,
\begin{equation}
    \label{error2D}
    E = h^2 \sum \limits_{i=1}^M \sum \limits_{j=1}^M \lvert u_{ij}^N - \bar u_{ij}^N \lvert \,.
\end{equation}
Moreover, the Experimental Order of Convergence (EOC) is computed using the errors from (\ref{error2D}) to check the expected order of the accuracy for the chosen examples. Note that to avoid a reduction of the accuracy due to boundary conditions when computing the EOC, we set exact values not only at the inflow boundary points, but also in a neighboring point outside of the computational interval to use the full stencil of the scheme for every inner grid point. 

\rv{
We also present the minimum or maximum values of numerical solutions at the final time $T$ as indicators of possible undershoots and overshoots in numerical solutions.
To assess the behavior of the fast sweeping method, we compute the $L_1$ norm ($||r||_1$) and the $L^{\infty}$ norm ($||r||_{\infty}$) of residuals for the algebraic systems of equations obtained with the considered numerical scheme.
}

Concerning the time steps, we choose maximal Courant numbers larger than allowed by a stability restriction of explicit schemes. In such a way, we document a stable behavior of the compact implicit schemes by preserving satisfactory accuracy. As the chosen examples do not contain any stiff features, if only the accuracy is considered, they can be solved better by explicit types of schemes. Nevertheless, our primary purpose is to show that no significant loss of accuracy is observed when applying large time steps that is in a large contrast with the unconditionally stable first order accurate implicit scheme.

The numerical methods and the graphical output are obtained using the Python programming language \cite{CS-R9526}.

\subsection{Linear advection equation}
First, we compute the linear advection equation given in (\ref{Linear2D}). As a guiding factor to compare the results, we compute the maximum value of Courant number in each direction,
\begin{equation}
    \label{courantnumber_advectionequation_experiments}
    C_{max}^x =  \frac{\tau}{h}\max \limits_{i,j} \{|v_{i+1/2,j}|\}\,, \hspace{1cm} C_{max}^y = \frac{\tau}{h}\max \limits_{i,j} \{|w_{i,j+1/2}|\} \,,
\end{equation}
over all edges of each finite volume. These "directional" Courant numbers will be larger than allowed by a stability restriction of explicit schemes. 


We choose the computational domain $x,y\in[-1, 1]$ and the velocity field representing the rotation in the form
\begin{equation}
    \label{velocityfield}
    \Vec{v} = (-2\pi y,2\pi x),
\end{equation}
where the exact solution for any initial function $u^0=u^0(x,y)$ is defined as
\begin{equation}
    \label{exactRot}
    u(x,y,t) = u^0\big(x \cos(2\pi t) + y \sin(2\pi t), y \cos(2\pi t) - x\sin(2\pi t)\big) \,.
\end{equation}

\subsubsection{First experiment}

In the first example, we consider the initial condition in a form of Gaussian,
\begin{equation}
    \label{Gausian2D}
    u^0(x,y) = e^{-10((x-0.25)^2+(y-0.25)^2)}
\end{equation}
and the final time $T=0.25$. In such a case, the Gaussian moves by a quarter of a cycle. The initial function and the exact solution are shown in Fig. \ref{FIG:firstexample_advection_gaussian}. For simplicity, the boundary conditions are set to 0.0 not only at the boundary points of the computational domain but also at neighboring points.

We perform experiments for three fixed choices of $\boldsymbol \omega$ to demonstrate the second order accuracy of the method. Using mesh sizes $M=40,80,160$ and $320$ with $N=M/5$ time steps (8, 16, 32, 64 sequentially), we maintain  $C_{max}^x = C_{max}^y = 3.92$.  The errors and the Experimental Orders of Convergence (EOCs) in Table \ref{TAB:firstexample_RotationGaussian} confirm the expected second order accuracy for the three pairs ${\omega}^{x,\pm}, {\omega}^{y,\pm} = 0, 1/2, 1$ used. 
We also apply the high-resolution scheme using ENO and WENO approximations to check the order of accuracy for such schemes. We also compute the errors and the EOCs using the first order accurate numerical scheme (\ref{firstorder}) to show the significant difference in the errors when using the higher order accurate schemes. 

\rv{
As shown in Table \ref{TAB:firstexample_RotationGaussian}, some minor oscillations in the numerical solution for certain constant values of $\boldsymbol \omega$ can be observed, especially for coarser meshes. These oscillations diminish with the mesh refinement and are eliminated when using the high-resolution schemes with ENO or WENO approximations.
We note that all presented errors and the minimum values are computed using only four sweeps \eqref{sweeps} to show that the expected EOC is obtained after this fixed number of GS iterations. 
}

\rv{
Additionally, we show the both norms of residuals after four iterations and, moreover, after another four iterations to show their behavior. Clearly, the norms are reduced significantly with more iterations for this example with smooth solution, especially for finer grids. Depending on preferences concerning the choice of the norm of residuals and its value at the end of iterations, any user can decide for a suitable stopping criteria of the iterative procedure. Note that comparisons between four and eight Gauss-Seidel iterations reveal no significant differences in errors, EOCs, or minimum values (beyond the fifth decimal place), so these results for eight iterations are omitted for this and next examples.
}

\ifshowfig
\begin{figure}[h]
\begin{center}
    \includegraphics[width=\textwidth]{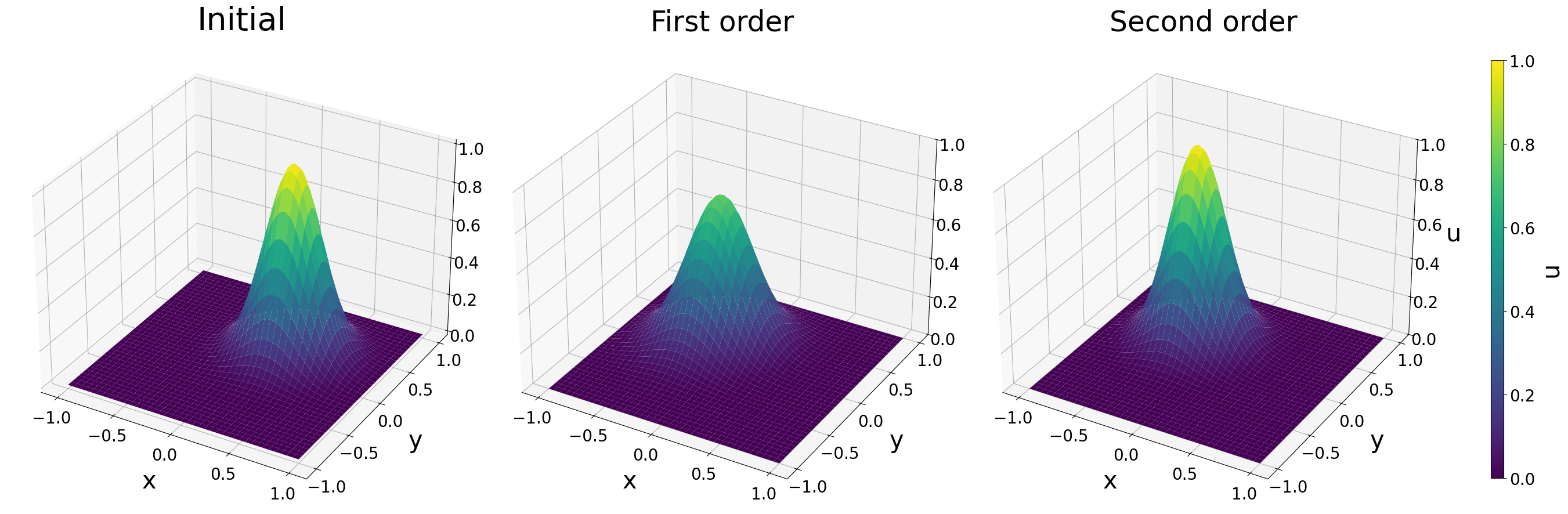}
    \end{center}
    \caption{The initial condition and the numerical solutions obtained by the first and second order accurate scheme with ${\omega}^{x,\pm}, {\omega}^{y,\pm} =1/2$ at $T=0.25$ for the coarsest mesh with $M=80$. }
    \label{FIG:firstexample_advection_gaussian}
\end{figure} 
\fi

\begin{table}[h]
    \begin{center}
    \begin{tabular}{c c | c c c c c }
    \multicolumn{2}{c}{} & \multicolumn{5}{c}{$1^{st}$ order, 4GS (8GS)} \\
    \hline

    $M $& $N$ & E      & EOC  & min                  & $||r||_1$         & $||r||_{\infty}$         \\ 
    \hline
    40  & 8   & 0.1727 & -    & 0.0 & $1.0 \cdot 10^{-3}$ ($6.3 \cdot 10^{-9}$) & $6.2 \cdot 10^{-3}$ ($1.1\cdot 10^{-7}$)\\
    80  & 16  & 0.1225 & 0.49 & 0.0 & $1.5 \cdot 10^{-4}$ ($1.0 \cdot 10^{-12}$)& $1.6 \cdot 10^{-3}$ ($1.6 \cdot 10^{-10}$)\\
    160 & 32  & 0.0766 & 0.67 & 0.0 & $1.3 \cdot 10^{-5}$ ($1.5 \cdot 10^{-16}$)& $3.4 \cdot 10^{-4}$ ($9.3 \cdot 10^{-14}$)\\
    320 & 64  & 0.0439 & 0.80 & 0.0 & $1.1 \cdot 10^{-6}$ ($5.0 \cdot 10^{-17}$)& $7.2 \cdot 10^{-5}$ ($4.4 \cdot 10^{-16}$)\\
    \hline \\
    \multicolumn{2}{c}{} & \multicolumn{5}{c}{$2^{nd}$ order, ${\omega}^{x,\pm}, {\omega}^{y,\pm} =0$, 4GS (8GS)} \\
    \hline
    $M $& $N$& E      & EOC  & min                  & $||r||_1$         & $||r||_{\infty}$ \\ 
    \hline
    40  & 8  & 0.1079 & -    & $-1.4 \cdot 10^{-2}$ & $5.9 \cdot 10^{-6}$ ($3.8 \cdot 10^{-12}$)& $6.3 \cdot 10^{-4}$ ($1.5 \cdot 10^{-10}$)  \\
    80  & 16  & 0.0355 & 1.60 & $-1.1 \cdot 10^{-4}$ & $4.6 \cdot 10^{-7}$ ($4.8 \cdot 10^{-16}$)& $1.2 \cdot 10^{-4}$ ($1.0 \cdot 10^{-13}$)  \\
    160 & 32  & 0.0096 & 1.88 & $-3.0 \cdot 10^{-5}$ & $4.1 \cdot 10^{-8}$ ($7.5 \cdot 10^{-17}$)& $2.5 \cdot 10^{-5}$ ($7.2 \cdot 10^{-16}$)  \\
    320 & 64 & 0.0024 & 1.97 & $-8.3 \cdot 10^{-6}$ & $4.4 \cdot 10^{-9}$ ($7.5 \cdot 10^{-17}$)& $5.4 \cdot 10^{-6}$ ($6.9 \cdot 10^{-16}$)  \\
    \hline \\
    \multicolumn{2}{c}{} & \multicolumn{5}{c}{$2^{nd}$ order, ${\omega}^{x,\pm}, {\omega}^{y,\pm} =1/2$, 4GS (8GS)} \\
    \hline
    $M $& $N$& E       & EOC  & min & $||r||_1$         & $||r||_{\infty}$\\ 
    \hline
    40  & 8  & 0.0415 & -    & $-1.8 \cdot 10^{-2}$ & $7.6 \cdot 10^{-5}$ ($3.2 \cdot 10^{-12}$)& $9.1 \cdot 10^{-4}$ ($1.6 \cdot 10^{-8}$)  \\
    80  & 16  & 0.0104 & 2.10 & $-5.6 \cdot 10^{-4}$ & $6.3 \cdot 10^{-6}$ ($2.0 \cdot 10^{-15}$)& $2.1 \cdot 10^{-5}$ ($5.0 \cdot 10^{-10}$)  \\
    160 & 32 & 0.0024 & 2.05 & 0.0 & $7.3 \cdot 10^{-7}$ ($7.2 \cdot 10^{-16}$)& $4.7 \cdot 10^{-5}$ ($6.3 \cdot 10^{-12}$)  \\
    320 & 64 & 0.0005 & 2.01 & 0.0 & $9.0 \cdot 10^{-8}$ ($1.0 \cdot 10^{-17}$)& $1.1 \cdot 10^{-6}$ ($4.8 \cdot 10^{-14}$)  \\
    \hline \\
    \multicolumn{2}{c}{} & \multicolumn{5}{c}{$2^{nd}$ order, ${\omega}^{x,\pm}, {\omega}^{y,\pm} =1$, 4GS (8GS)} \\
    \hline
    $M $& $N$& E       & EOC  & min & $||r||_1$         & $||r||_{\infty}$  \\ 
    \hline
    40  & 8  & 0.0742 & -    & $-2.8 \cdot 10^{-2}$  & $4.4 \cdot 10^{-4}$ ($1.4 \cdot 10^{-8}$) & $5.0 \cdot 10^{-3}$ ($5.4 \cdot 10^{-7}$)  \\
    80  & 16  & 0.0252 & 1.55 & $-6.4 \cdot 10^{-3}$  & $2.0 \cdot 10^{-5}$ ($1.1 \cdot 10^{-10}$)& $8.9 \cdot 10^{-4}$ ($1.9 \cdot 10^{-8}$)  \\
    160 & 32 & 0.0067 & 1.91 & $-2.5 \cdot 10^{-5}$  & $2.6 \cdot 10^{-6}$ ($4.5 \cdot 10^{-13}$)& $2.0 \cdot 10^{-4}$ ($3.0 \cdot 10^{-10}$)  \\
    320 & 64 & 0.0017 & 1.98 & 0.0 & $3.5 \cdot 10^{-7}$ ($1.1 \cdot 10^{-15}$)& $5.1 \cdot 10^{-5}$ ($2.6 \cdot 10^{-12}$)  \\
    \hline \\
    \multicolumn{2}{c}{} & \multicolumn{5}{c}{ENO, 4GS (8GS)} \\
    \hline
    $M $& $N$& E       & EOC  & min & $||r||_1$         & $||r||_{\infty}$ \\ 
    \hline
    40  & 8  & 0.0753 & -    & 0.0 & $3.3 \cdot 10^{-5}$ ($6.7 \cdot 10^{-11}$)& $3.6 \cdot 10^{-4}$ ($2.7 \cdot 10^{-9}$)  \\
    80  & 16  & 0.0272 & 1.46 & 0.0 & $1.4 \cdot 10^{-6}$ ($1.6 \cdot 10^{-12}$)& $4.0 \cdot 10^{-5}$ ($9.0 \cdot 10^{-11}$)  \\
    160 & 32 & 0.0086 & 1.65 & 0.0 & $3.3 \cdot 10^{-8}$ ($1.7 \cdot 10^{-14}$)& $1.8 \cdot 10^{-6}$ ($9.6 \cdot 10^{-12}$)  \\
    320 & 64 & 0.0025 & 1.76 & 0.0 & $1.1 \cdot 10^{-9}$ ($3.0 \cdot 10^{-14}$)& $3.4 \cdot 10^{-7}$ ($4.7 \cdot 10^{-11}$)  \\
    \hline \\
    \multicolumn{2}{c}{} & \multicolumn{5}{c}{WENO, 4GS (8GS)} \\
    $M $& $N$& E       & EOC  & min  & $||r||_1$         & $||r||_{\infty}$  \\ 
    \hline
    40  & 8  & 0.0673 & -    & 0.0 & $9.6 \cdot 10^{-6}$ ($2.9 \cdot 10^{-10}$)& $1.2 \cdot 10^{-4}$ ($3.3 \cdot 10^{-9}$)  \\
    80  & 16 & 0.0217 & 1.63 & 0.0 & $4.5 \cdot 10^{-7}$ ($2.2 \cdot 10^{-11}$)& $1.1 \cdot 10^{-5}$ ($2.6 \cdot 10^{-10}$)  \\
    160 & 32 & 0.0061 & 1.82 & 0.0 & $2.4 \cdot 10^{-8}$ ($9.0 \cdot 10^{-14}$)& $3.7 \cdot 10^{-6}$ ($6.6 \cdot 10^{-12}$)  \\
    320 & 64 & 0.0015 & 2.01 & 0.0 & $3.6 \cdot 10^{-9}$ ($2.3 \cdot 10^{-14}$)& $3.5 \cdot 10^{-7}$ ($7.8 \cdot 10^{-12}$)  \\
    \hline 
    \end{tabular}
    \end{center}
    \caption{The numerical errors (E), EOCs, minimum values, the norm of residuals $||r||_1$ and $||r||_{\infty}$ of the first and second order schemes and the high-resolution scheme with ENO and WENO approximations for the rotation of Gaussian with $C_{max}^x = C_{max}^y = 3.92$. The norms of residuals are shown for four (4GS) and eight (8GS) Gauss-Seidel iterations. }
    \label{TAB:firstexample_RotationGaussian}
\end{table}

\clearpage
\subsubsection{Second experiment}
Next, we present numerical results for a non-smooth solution. As an initial condition for $x,y\in [-1,1]$, we choose four formations in four quadrants, see (\ref{InitialSpec}) and Fig. \ref{FIG:discSpecInitialExact}. In particular, a Gaussian with the center at $(0.5,0.5)$ and the height equals $1$; a cone with the center at $(-0.5,0.5)$ with the maximal radius $r_{max}=0.25$ and the height 1; a half sphere with the center at $(-0.5,-0.5)$ and the maximal radius $r_{max}$ and the height 1; and a circle with a center at $(0.5, -0.5)$ with radius $r_{max}$ with the value 1 inside the circle and the value 0 everywhere else. 
\begin{align}  
    \label{InitialSpec}
    u^0(x, y) =
    \begin{cases}
    e^{100 \left(-(x - 0.5)^2 - (y - 0.5)^2\right) } & \text{if } x \geq 0 \text{ and } y \geq 0 \text{ and } (x - 0.5)^2 + (y - 0.5)^2 < 0.3^2, \\
    1 - \frac{\sqrt{(x + 0.5)^2 + (y - 0.5)^2}}{0.25} & \text{if } x < 0 \text{ and } y \geq 0 \text{ and } \sqrt{(x + 0.5)^2 + (y - 0.5)^2} \leq 0.25, \\
    \sqrt{1 - \Big(\frac{\sqrt{(x + 0.5)^2 + (y + 0.5)^2}}{0.25} \Big)^2} & \text{if } x < 0 \text{ and } y < 0 \text{ and } \sqrt{(x + 0.5)^2 + (y + 0.5)^2} \leq 0.25, \\
    1 & \text{if } x \geq 0 \text{ and } y < 0 \text{ and } \sqrt{(x - 0.5)^2 + (y + 0.5)^2} \leq 0.25, \\
    0 & \text{otherwise}.
    \end{cases}
\end{align}
Note that the maximum value equals 1 and the minimum equals 0.

The velocity field is the same as in (\ref{velocityfield}), and we choose the final time $T=0.25$. The exact solution is obtained using (\ref{exactRot}). The initial condition and the final solution are shown in \mbox{Figure \ref{FIG:discSpecInitialExact}}. We choose $M=40,80,160,320$ and $N=M/$ time steps (specifically 8, 16, 32 and 64), leading to the same Courant number as in the first example, $C_{max}^x = C_{max}^y = 3.92$.

This time, it is necessary to use the high-resolution scheme with ENO and WENO approximations to suppress oscillations in numerical solutions. The errors and EOCs are shown in Table \ref{TAB:secondexample_advection_rotation}, together with the comparison to the errors when using the first order accurate scheme. 

\rv{
The numerical results are presented in the analogous form as in the previous example. All errors and the minimum and maximum values are presented for numerical solutions obtained with four sweeps \eqref{sweeps}, but the norms of residuals are additionally given also after another four GS iterations. From the Figure \ref{FIG:discSpecInitialExact} and the Table \ref{TAB:secondexample_advection_rotation} one can clearly see improvements in accuracy, especially with respect to the first order scheme, with the most accurate results obtained using the WENO approximations.
}

\ifshowfig
\begin{figure}[h]
    \begin{center} 
    \includegraphics[width=\textwidth]{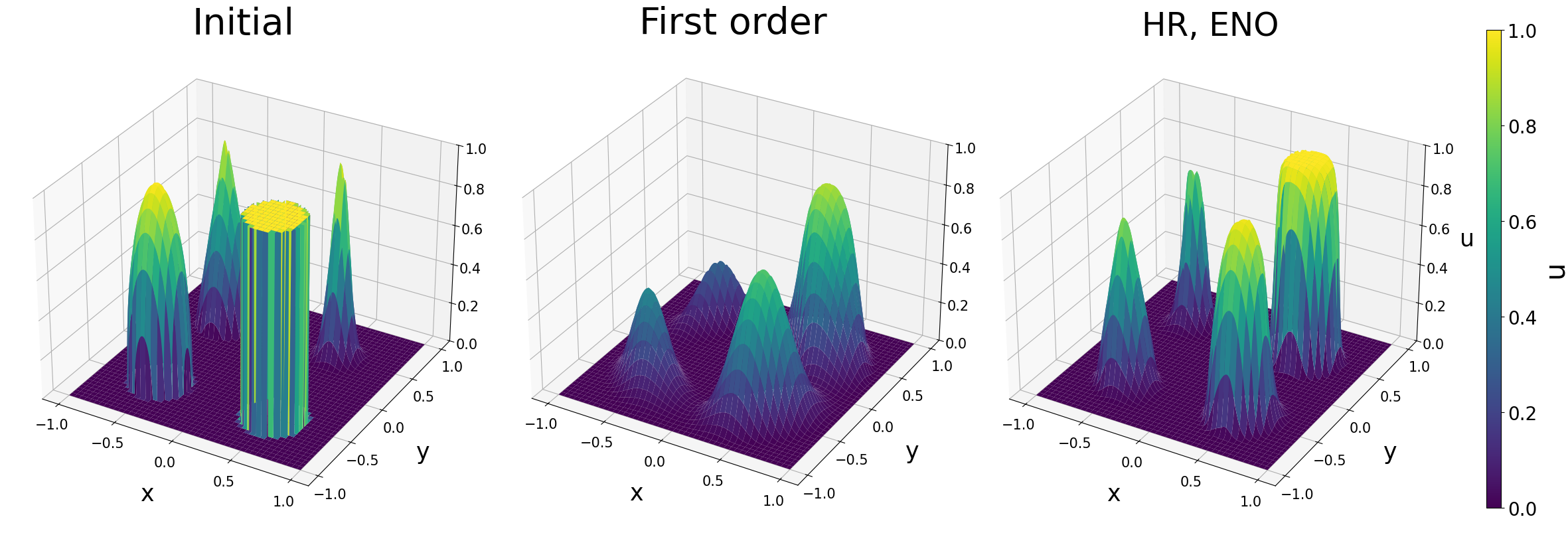}
    \end{center}
    \caption{The initial condition for the example of the rotation of \eqref{InitialSpec} and numerical results using the first order accurate scheme and high-resolution scheme with ENO approximations for $C_{max}^x = C_{max}^y = 3.92$ on the finest mesh with $M=320$.}
    \label{FIG:discSpecInitialExact}
\end{figure} 
\fi 

\begin{table}[h]
    \begin{center}
    \begin{tabular}{c c | c c c c c c}
    \multicolumn{2}{c}{} & \multicolumn{6}{c}{$1^{st}$ order, 4GS (8GS)} \\
    \hline
    $M $& $N$ & E      & EOC  & min     & max           & $||r||_1$         & $||r||_{\infty}$          \\ 
    \hline
    40  & 8   & 0.5610 & -    & 0.0 & 0.32 & $8.1 \cdot 10^{-4}$ ($3.1 \cdot 10^{-10}$)& $7.7 \cdot 10^{-3}$ ($1.7 \cdot 10^{-9}$)\\
    80  & 16  & 0.5421 & 0.04 & 0.0 & 0.40 & $1.3 \cdot 10^{-4}$ ($3.2 \cdot 10^{-16}$)& $1.6 \cdot 10^{-3}$ ($2.9 \cdot 10^{-14}$)\\
    160 & 32  & 0.4819 & 0.17 & 0.0 & 0.54 & $8.5 \cdot 10^{-6}$ ($1.0 \cdot 10^{-16}$)& $1.7 \cdot 10^{-4}$ ($4.4 \cdot 10^{-16}$)\\
    320 & 64  & 0.3855 & 0.32 & 0.0 & 0.72 & $2.5 \cdot 10^{-7}$ ($1.1 \cdot 10^{-16}$)& $8.9 \cdot 10^{-6}$ ($8.8 \cdot 10^{-16}$)\\
    \hline \\
    \multicolumn{2}{c}{} & \multicolumn{6}{c}{ENO, 4GS (8GS)} \\
    \hline
    $M $& $N$& E       & EOC  & min  & max & $||r||_1$         & $||r||_{\infty}$ \\ 
    \hline
    40  & 2  & 0.5186 & -    & 0.0                   & 0.49 & $2.0 \cdot 10^{-3}$ ($1.9 \cdot 10^{-7}$)& $7.2 \cdot 10^{-2}$ ($1.6 \cdot 10^{-6}$)  \\
    80  & 4  & 0.3755 & 0.46 & 0.0                   & 0.72 & $2.4 \cdot 10^{-5}$ ($1.2 \cdot 10^{-9}$)& $9.2 \cdot 10^{-4}$ ($8.1 \cdot 10^{-8}$)  \\
    160 & 8  & 0.2231 & 0.75 & 0.0                   & 0.93 & $3.7 \cdot 10^{-6}$ ($1.3 \cdot 10^{-11}$)& $5.0 \cdot 10^{-4}$ ($2.5 \cdot 10^{-9}$)  \\
    320 & 16 & 0.1228 & 0.86 & $-1.7 \cdot 10^{-17}$ & 0.99 & $1.1 \cdot 10^{-8}$ ($7.7 \cdot 10^{-14}$)& $4.4 \cdot 10^{-6}$ ($4.6 \cdot 10^{-11}$)  \\
    \hline \\
    \multicolumn{2}{c}{} & \multicolumn{6}{c}{WENO, 4GS (8GS)} \\
    \hline
    $M $& $N$& E       & EOC  & min  & max & $||r||_1$         & $||r||_{\infty}$ \\ 
    \hline
    40  & 2  & 0.4991 & -    & 0.0                   & 0.55 & $1.9 \cdot 10^{-4}$ ($1.1 \cdot 10^{-8}$)& $2.2 \cdot 10^{-3}$ ($1.3 \cdot 10^{-7}$)  \\
    80  & 4  & 0.3485 & 0.51 & 0.0                   & 0.76 & $1.4 \cdot 10^{-5}$ ($9.9 \cdot 10^{-10}$)& $6.9 \cdot 10^{-4}$ ($3.6 \cdot 10^{-8}$)  \\
    160 & 8  & 0.1949 & 0.83 & $-4.0 \cdot 10^{-10}$ & 0.95 & $2.3 \cdot 10^{-6}$ ($6.7 \cdot 10^{-11}$)& $4.7 \cdot 10^{-5}$ ($1.3 \cdot 10^{-9}$)  \\
    320 & 16 & 0.1016 & 0.94 & $-1.6 \cdot 10^{-8} $ & 0.99 & $4.0 \cdot 10^{-7}$ ($1.2 \cdot 10^{-11}$)& $5.6 \cdot 10^{-5}$ ($1.3 \cdot 10^{-9}$)  \\
    \hline 
    \end{tabular}
    \end{center}
    \caption{The numerical errors (E), the EOCs, minimum and maximum values, the norm of residuals $||r||_1$ and $||r||_{\infty}$  of the first order scheme and the high-resolution scheme with ENO and WENO approximations for the rotation of four shapes \eqref{InitialSpec} with $C_{max}^x = C_{max}^y = 3.92$. The norms of residuals are shown for four (4GS) and eight (8GS) Gauss-Seidel iterations. }
    \label{TAB:secondexample_advection_rotation}
\end{table}

\clearpage

\subsection{Burgers' equation}
We compute a representative nonlinear problem in the form of the Burgers equation,
\begin{equation}
    \label{Burgers2D}
    \partial_{t} u +  \partial_x \Big(\frac{u^2}{2}\Big) +  \partial_y \Big(\frac{u^2}{2}\Big) = 0 \,.
\end{equation}
In this case, the maximum Courant number is defined as $C_{max}=\frac{\tau}{h}u_{max}$, with $u_{max}=\max(|u_{i,j}^n|)$ over all time steps and finite volumes.

\subsubsection{First example}
For the first 2D example we choose the computation domain $x,y\in[-1, 1]$ and $t\in[0, 0.5]$ with a smooth initial condition (Fig. \ref{FIG:2D_nonlinearsinus}) in the form
\begin{equation}
    \label{sinus_nonlinear_initial}
    u^0(x,y) = \sin(\pi x)  \sin(\pi y) / 2 \,.
\end{equation} 
The exact solution can be obtained numerically using the method of characteristics by solving the algebraic equations for $u=u(x_i,y_j,t^n)$
\begin{equation}
    \label{numericalsolution_sinus_nonlinear}
     u = \sin(\pi (x_i - ut^n)) \sin(\pi (y_j - ut^n)) / 2 \,.
\end{equation}
For this example, we chose sequentially $M=80,160,320$ and $640$, and the number of time steps is $N=M/40$ corresponding to $N=$. Together with the maximum absolute value (and also the minimum) of the function $u^0$ being 0.5, one obtains $C_{max} = 5$.

We present results obtained using the compact implicit numerical scheme with constant values of ${\omega}^{x,\pm}, {\omega}^{y,\pm} = 0, 1/2, 1$, respectively. The errors \eqref{error2D} at the final time $T = 0.5$, the EOCs, and the minimum and maximum values of $u$ (computed with four sweeps) are presented at the final time T, revealing slight undershoots and overshoots. The errors and EOCs in Table \ref{TAB:firstexample_burgers_sinus} confirm the expected order of accuracy for the different choices of ${\omega}^{x,\pm}, {\omega}^{y,\pm}$.

We also present results using the high-resolution scheme with ENO and WENO approximations to demonstrate the order of accuracy. Notably, the maximum and minimum values of $u$ obtained with the ENO and WENO approximations remain within the interval $(-0.5, 0.5)$, in contrast to the results with constant values of $\boldsymbol{\omega}$ (see Table \ref{TAB:firstexample_burgers_sinus}).

\rv{
The residual norms, $||r||_1$ and $||r||_{\infty}$, are also monitored, both of which decrease with mesh refinement and with additional GS iterations. Comparisons between four and eight Gauss-Seidel iterations reveal no significant differences in errors, EOCs, or minimum and maximum values, as observed before for previous examples. 
These results may be compared also with the errors of the first order accurate scheme presented in Table \ref{TAB:firstexample_burgers_sinus}.
}

\ifshowfig
\begin{figure}[h]
    \begin{center}
    \includegraphics[width=\textwidth]{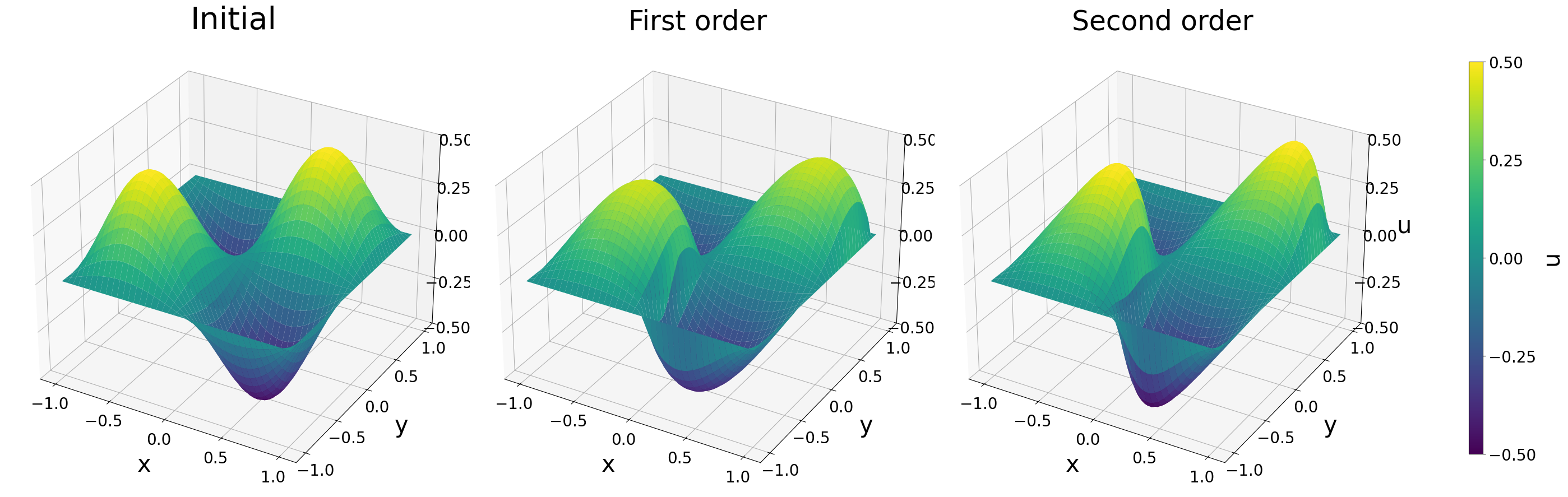}
    \end{center}
    \caption{The initial condition (left), and the numerical solution obtained by the first (middle) and the second (right) order scheme at $T=0.5$ for the nonlinear problem with the smooth initial condition \eqref{sinus_nonlinear_initial}, ${\omega}^{x,\pm}, {\omega}^{y,\pm} = 1/2$, $C_{max} = 5$, and $M=80$ (the coarsest mesh). }
    \label{FIG:2D_nonlinearsinus}
\end{figure} 
\fi

\begin{table}[h]
    \begin{center}
    \begin{tabular}{c c | c c c c c c }
    \multicolumn{2}{c}{} & \multicolumn{6}{c}{$1^{st}$ order, 4GS (8GS)} \\
    \hline
    $M$ & $N$ & E      & EOC  & min    & max   & $||r||_1$                                & $||r||_{\infty}$          \\ 
    \hline
    80  & 2   & 0.0944 & -    & -0.427 & 0.427 & $4.4 \cdot 10^{-5}$ ($9.6 \cdot 10^{-7}$)& $8.6 \cdot 10^{-3}$ ($2.7\cdot 10^{-4}$)\\
    160 & 4   & 0.0618 & 0.61 & -0.449 & 0.449 & $7.7 \cdot 10^{-6}$ ($5.3 \cdot 10^{-8}$)& $4.7 \cdot 10^{-3}$ ($9.4 \cdot 10^{-5}$)\\
    320 & 8   & 0.0378 & 0.70 & -0.467 & 0.467 & $7.5 \cdot 10^{-7}$ ($1.0 \cdot 10^{-9}$)& $1.0 \cdot 10^{-3}$ ($7.8 \cdot 10^{-6}$)\\
    640 & 16  & 0.0214 & 0.82 & -0.481 & 0.481 & $4.3 \cdot 10^{-9}$ ($2.5 \cdot 10^{-16}$)& $1.6 \cdot 10^{-5}$ ($7.7 \cdot 10^{-16}$)\\
    \hline \\
    \multicolumn{2}{c}{} & \multicolumn{6}{c}{$2^{nd}$ order, ${\omega}^{x,\pm}, {\omega}^{y,\pm} =0$, 4GS (8GS)} \\
    \hline
    $M $& $N$& E      & EOC  & min    & max   & $||r||_1$                                  & $||r||_{\infty}$  \\ 
    \hline
    80  & 2  & 0.0280 & -    & -0.511 & 0.511 & $1.5 \cdot 10^{-8}$ ($2.9 \cdot 10^{-16}$) & $4.7\cdot 10^{-6}$ ($9.9 \cdot 10^{-16}$)  \\
    160 & 4  & 0.0097 & 1.52 & -0.506 & 0.506 & $2.6 \cdot 10^{-10}$ ($2.8 \cdot 10^{-16}$)& $1.9 \cdot 10^{-7}$ ($1.2 \cdot 10^{-15}$)  \\
    320 & 8  & 0.0029 & 1.74 & -0.502 & 0.502 & $8.9 \cdot 10^{-13}$ ($2.8 \cdot 10^{-16}$)& $1.6 \cdot 10^{-9}$ ($1.3 \cdot 10^{-15}$)  \\
    640 & 16 & 0.0007 & 1.88 & -0.500 & 0.500 & $1.3 \cdot 10^{-15}$ ($2.8 \cdot 10^{-16}$)& $4.1 \cdot 10^{-12}$ ($1.3 \cdot 10^{-15}$)  \\
    \hline \\
    \multicolumn{2}{c}{} & \multicolumn{6}{c}{$2^{nd}$ order, ${\omega}^{x,\pm}, {\omega}^{y,\pm} =1/2$, 4GS (8GS)} \\
    \hline
    $M $& $N$& E       & EOC  & min  & max & $||r||_1$                                         & $||r||_{\infty}$  \\ 
    \hline
    80  & 2  & 0.0208 & -    & -0.507 & 0.507 & $1.3 \cdot 10^{-7}$ ($6.6 \cdot 10^{-12}$) & $9.3 \cdot 10^{-6}$ ($3.2 \cdot 10^{-10}$)  \\
    160 & 4  & 0.0070 & 1.56 & -0.504 & 0.504 & $1.0 \cdot 10^{-9}$ ($6.0 \cdot 10^{-15}$) & $8.7 \cdot 10^{-8}$ ($5.3 \cdot 10^{-13}$)  \\
    320 & 8  & 0.0020 & 1.79 & -0.501 & 0.501 & $8.0 \cdot 10^{-12}$ ($3.2 \cdot 10^{-16}$)& $1.1 \cdot 10^{-9}$ ($1.6 \cdot 10^{-15}$)  \\
    640 & 16 & 0.0005 & 1.92 & -0.500 & 0.500 & $6.8 \cdot 10^{-14}$ ($3.1 \cdot 10^{-16}$)& $1.9 \cdot 10^{-11}$ ($1.4 \cdot 10^{-15}$)  \\
    \hline \\
    \multicolumn{2}{c}{} & \multicolumn{6}{c}{$2^{nd}$ order, ${\omega}^{x,\pm}, {\omega}^{y,\pm} =1$, 4GS (8GS)} \\
    \hline
    $M $& $N$& E       & EOC  & min  & max & $||r||_1$           & $||r||_{\infty}$  \\ 
    \hline
    80  & 2  & 0.0148 & -    & -0.503 & 0.503 & $6.4 \cdot 10^{-7}$ ($1.2 \cdot 10^{-10}$) & $1.6\cdot 10^{-5}$ ($5.3 \cdot 10^{-9}$)  \\
    160 & 4  & 0.0050 & 1.56 & -0.502 & 0.502 & $6.2 \cdot 10^{-9}$ ($1.5 \cdot 10^{-13}$) & $3.5 \cdot 10^{-7}$ ($1.2 \cdot 10^{-11}$)  \\
    320 & 8  & 0.0014 & 1.80 & -0.500 & 0.500 & $6.5 \cdot 10^{-11}$ ($6.2 \cdot 10^{-16}$)& $1.7 \cdot 10^{-8}$ ($2.4 \cdot 10^{-14}$)  \\
    640 & 16 & 0.0003 & 1.94 & -0.500 & 0.500 & $5.1 \cdot 10^{-13}$ ($4.9 \cdot 10^{-16}$)& $1.4 \cdot 10^{-10}$ ($2.4 \cdot 10^{-15}$)  \\
    \hline \\
    \multicolumn{2}{c}{} & \multicolumn{6}{c}{ENO, 4GS (8GS)} \\
    \hline
    $M $& $N$& E       & EOC  & min  & max & $||r||_1$         & $||r||_{\infty}$  \\ 
    \hline
    80  & 2  & 0.0267 & -    & -0.484 & 0.484 & $1.1 \cdot 10^{-5}$ ($4.5 \cdot 10^{-8}$)& $2.7 \cdot 10^{-4}$ ($2.1 \cdot 10^{-6}$)  \\
    160 & 4  & 0.0087 & 1.60 & -0.494 & 0.494 & $7.1 \cdot 10^{-8}$ ($9.8 \cdot 10^{-11}$)& $4.5 \cdot 10^{-6}$ ($1.6 \cdot 10^{-8}$)  \\
    320 & 8  & 0.0030 & 1.52 & -0.498 & 0.498 & $6.9 \cdot 10^{-11}$ ($6.0 \cdot 10^{-16}$)& $1.2 \cdot 10^{-7}$ ($9.8 \cdot 10^{-14}$)  \\
    640 & 16 & 0.0008 & 1.82 & -0.499 & 0.499 & $6.3 \cdot 10^{-13}$ ($3.9 \cdot 10^{-16}$)& $1.9 \cdot 10^{-9}$ ($2.9 \cdot 10^{-15}$)  \\
    \hline \\
    \multicolumn{2}{c}{} & \multicolumn{6}{c}{WENO, 4GS (8GS)} \\
    $M $& $N$& E       & EOC  & min  & max & $||r||_1$           & $||r||_{\infty}$  \\ 
    \hline
    80  & 2  & 0.0252 & -    & -0.485 & 0.485 & $1.0 \cdot 10^{-5}$  ($4.4 \cdot 10^{-8}$) & $2.0 \cdot 10^{-4}$ ($1.8 \cdot 10^{-6}$)  \\
    160 & 4  & 0.0080 & 1.64 & -0.495 & 0.495 & $7.6 \cdot 10^{-8}$  ($4.0 \cdot 10^{-11}$)& $8.5 \cdot 10^{-6}$ ($5.1 \cdot 10^{-9}$)  \\
    320 & 8  & 0.0029 & 1.43 & -0.498 & 0.498 & $2.4 \cdot 10^{-11}$ ($3.8 \cdot 10^{-16}$)& $4.3 \cdot 10^{-8}$ ($1.0 \cdot 10^{-14}$)  \\
    640 & 16 & 0.0008 & 1.83 & -0.499 & 0.499 & $4.0 \cdot 10^{-13}$ ($3.1 \cdot 10^{-16}$)& $1.0 \cdot 10^{-8}$ ($3.9 \cdot 10^{-15}$)  \\
    \hline 
    \end{tabular}
    \end{center}
    \caption{The numerical errors (E), EOCs, minimum values, the norm of residuals $||r||_1$ and $||r||_{\infty}$  of the first and second order schemes and the high-resolution scheme with ENO and WENO approximations for the Burgers equation with the initial sine function with $C_{max}^x = C_{max}^y = 5.0$. The norms of residuals are shown for four (4GS) and eight (8GS) Gauss-Seidel iterations. }
    \label{TAB:firstexample_burgers_sinus}
\end{table}

\clearpage
\subsubsection{Second example}
The second problem is, again, defined for $x,y\in[-1,1]$ and will consider a rarefaction wave effect. The initial condition, see Figure \ref{FIG:exact_rarefraction}, is given in the form 
\begin{align}
    \label{initial_rarefraction}
     u^0(x,y) = \begin{cases}
         u_L \hspace{0.5cm} &\text{ if } \frac{x + y}{2}<0 \\
        u_R  \hspace{0.5cm} &\text{ otherwise }
     \end{cases} \,,
\end{align} 
for $u_L = -1$ and $u_R = 1$; and, the exact solution is prescribed as
\begin{align}
    \label{exact_rarefraction}
     u(x, y, t) =
        \begin{cases}
        u_L \hspace{0.5cm}& \text{if } \frac{x + y}{2} \leq - t \\
        \frac{x + y}{2t} \hspace{0.5cm}& \text{if } - t < \frac{x + y}{2} \leq  t \\
        u_R \hspace{0.5cm}& \text{if } \frac{x + y}{2} >t
        \end{cases} \,,
\end{align}
and we compute the solution for $t\in[0,0.4]$. 

\ifshowfig
\begin{figure}[h]
    \begin{center}
    \includegraphics[width=\textwidth]{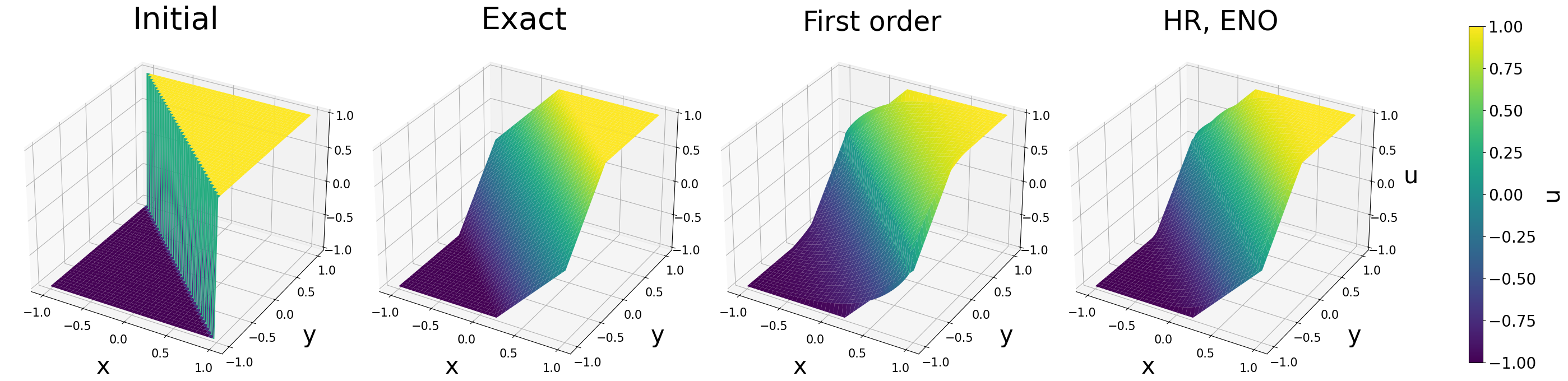}
    \end{center}
    \caption{The initial condition (left), the exact solution (second to the left) and the numerical solution using the first order accurate scheme and the high-resolution scheme with ENO approximation (third and fourth figure, respectively) at $T=0.4$  with the coarsest mesh $M=80$ for the nonlinear problem with the initial condition (\ref{exact_rarefraction}). }
    \label{FIG:exact_rarefraction}
\end{figure} 
\fi 

This time we choose $N=M/20$ number of time steps for $M = 40, 80, 160, 320$ (2, 4, 8, 16 corresponding numbers of time steps), with the value $C_{max} = 4.0$. The high-resolution scheme with ENO and WENO (with $\bar{\omega}=1/3$) approximation were used to compute the solution, together with the first order accurate scheme, to demonstrate the difference in the errors obtained for four Gauss-Seidel iterations, see Table \ref{TAB:secondexample_burgers_rarefaction}.

\begin{table}[h]
    \begin{center}
    \begin{tabular}{c c | c c c c }
    \multicolumn{2}{c}{} & \multicolumn{4}{c}{$1^{st}$ order, 4GS (8GS) [12GS]} \\
    \hline
    $M $& $N$ & E     & EOC  & $||r||_1$         & $||r||_{\infty}$         \\ 
    \hline
    40  & 2  & 0.5035 & -    & $7.0 \cdot 10^{-16}$ ($6.9 \cdot 10^{-16}$) [$6.9 \cdot 10^{-16}$] & $1.0 \cdot 10^{-15}$ ($1.0 \cdot 10^{-15}$) [$9.9 \cdot 10^{-16}$]\\
    80  & 4  & 0.3762 & 0.37 & $6.8 \cdot 10^{-16}$ ($6.8 \cdot 10^{-16}$) [$6.9 \cdot 10^{-16}$] & $9.9 \cdot 10^{-16}$ ($9.9 \cdot 10^{-16}$) [$9.9 \cdot 10^{-16}$]\\
    160 & 8  & 0.2613 & 0.52 & $4.8 \cdot 10^{-9}$  ($7.0 \cdot 10^{-16}$) [$6.9 \cdot 10^{-16}$] & $9.1 \cdot 10^{-6}$  ($1.0 \cdot 10^{-15}$) [$1.0 \cdot 10^{-15}$]\\
    320 & 16 & 0.1732 & 0.59 & $2.0 \cdot 10^{-10}$ ($6.9 \cdot 10^{-16}$) [$6.9 \cdot 10^{-16}$] & $1.4 \cdot 10^{-7}$  ($1.2 \cdot 10^{-15}$) [$1.2 \cdot 10^{-15}$]\\
    \hline \\
    \multicolumn{2}{c}{} & \multicolumn{4}{c}{ENO, 4GS (8GS) [12GS]} \\
    \hline
    $M $& $N$& E     & EOC  & $||r||_1$         & $||r||_{\infty}$\\ 
    \hline
    40  & 2 & 0.2283 & -    & $1.3 \cdot 10^{-3}$  ($6.0 \cdot 10^{-9}$)  [$1.5 \cdot 10^{-13}$] & $7.1 \cdot 10^{-3}$ ($3.2 \cdot 10^{-8}$)  [$1.5 \cdot 10^{-13}$] \\
    80  & 4 & 0.1362 & 0.74 & $5.0 \cdot 10^{-5}$  ($4.9 \cdot 10^{-8}$)  [$7.2 \cdot 10^{-10}$] & $5.0 \cdot 10^{-4}$ ($7.1 \cdot 10^{-7}$)  [$7.2 \cdot 10^{-10}$] \\
    160 & 8 & 0.0751 & 0.85 & $2.5 \cdot 10^{-6}$  ($7.6 \cdot 10^{-10}$) [$2.6 \cdot 10^{-11}$] & $5.7 \cdot 10^{-5}$ ($2.6 \cdot 10^{-8}$)  [$2.6 \cdot 10^{-11}$] \\
    320 & 16& 0.0396 & 0.92 & $9.5 \cdot 10^{-10}$ ($2.5 \cdot 10^{-14}$) [$6.7 \cdot 10^{-15}$] & $7.4 \cdot 10^{-7}$ ($6.0 \cdot 10^{-11}$) [$6.9 \cdot 10^{-15}$] \\
    \hline \\
    \multicolumn{2}{c}{} & \multicolumn{4}{c}{WENO, 4GS (8GS) [12GS]} \\
    \hline
    $M $& $N$& E     & EOC  & $||r||_1$         & $||r||_{\infty}$  \\ 
    \hline
    40  & 2 & 0.2087 & -    & $1.4 \cdot 10^{-3}$ ($4.1 \cdot 10^{-8}$)  [$2.0 \cdot 10^{-13}$] & $7.8 \cdot 10^{-3}$ ($2.1 \cdot 10^{-7}$)  [$9.4 \cdot 10^{-13}$] \\
    80  & 4 & 0.1225 & 0.76 & $6.8 \cdot 10^{-5}$ ($7.5 \cdot 10^{-8}$)  [$5.8 \cdot 10^{-11}$] & $5.4 \cdot 10^{-4}$ ($8.9 \cdot 10^{-7}$)  [$5.3 \cdot 10^{-10}$] \\
    160 & 8 & 0.0670 & 0.86 & $1.6 \cdot 10^{-6}$ ($4.6 \cdot 10^{-10}$) [$4.2 \cdot 10^{-13}$] & $5.1 \cdot 10^{-5}$ ($1.4 \cdot 10^{-8}$)  [$1.4 \cdot 10^{-11}$] \\
    320 & 16& 0.0350 & 0.93 & $8.1 \cdot 10^{-9}$ ($4.3 \cdot 10^{-13}$) [$8.4 \cdot 10^{-16}$] & $3.0 \cdot 10^{-6}$ ($3.1 \cdot 10^{-10}$) [$3.9 \cdot 10^{-14}$] \\
    \hline 
    \end{tabular}
    \end{center}
    \caption{The numerical errors (E), the EOCs, minimum and maximum values, the norm of residuals $||r||_1$ and $||r||_{\infty}$ of the first order scheme and the high-resolution scheme with ENO and WENO approximations for the rarefactrion wave with $C_{max}^x = C_{max}^y = 4.0$. The results are shown for four (4GS), eight (8GS) and twelve (12GS) Gauss-Seidel iterations. }
    \label{TAB:secondexample_burgers_rarefaction}
\end{table}

In this case, the maximum and minimum values are omitted from Table \ref{TAB:secondexample_burgers_rarefaction} since no oscillations were observed (the solution remains bounded within [-1,1]).

\rv{
The residual norms, $||r||_1$ and $||r||_{\infty}$ are monitored in Table \ref{TAB:secondexample_burgers_rarefaction} and shown to decrease with mesh refinement and additional iterations. As now the problem with the sonic point \cite{shu_essentially_1998,leveque_finite_2004} become more challenging, we give comparisons between four, eight, and twelve Gauss-Seidel iterations. The convergence rate of the residual norms is slower, but these norms can be clearly significantly reduced, especially for finer meshes.
}

These results can be compared with those of the first order accurate scheme presented in Table \ref{TAB:secondexample_burgers_rarefaction}.

\clearpage
\subsubsection{Third example}
\label{example-burgers-shocks}
\begin{figure}[h]
    \begin{center}
    \includegraphics[width=\textwidth]{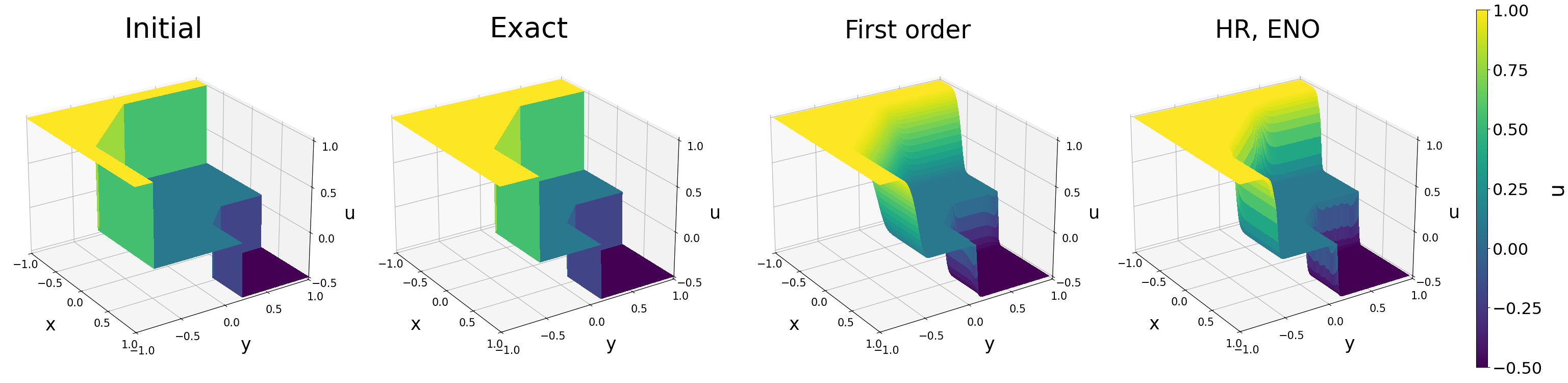}
    \end{center}
    \caption{The visualization of the initial condition (left), the exact solution (second left) and the the numerical solutions (right) obtained by the first order accurate and high-resolution scheme with ENO approximation at $T=0.4$ for the coarser mesh with $M=160$ for the nonlinear problem with the initial condition \eqref{2D-initial_shock}. }
    \label{FIG:2D-exact_shock}
\end{figure}

The next problem deals with a discontinuity in the form of a shock wave and is defined for $x,y\in[-1,1]$ with the computational time $t\in[0,0.4]$. The initial condition in Figure \ref{FIG:2D-exact_shock} is given in the form 
\begin{align}
    \label{2D-initial_shock}
     u^0(x,y) = \begin{cases}
        u_{L1} & \text{if } x < -0.8 \text{ or } y < -0.8 \text{ or } x + y < -0.8 \\
        u_{L2} & \text{if } x <  0.2 \text{ or } y <  0.2 \text{ or } x + y <  0.7 \\
        u_{R2} & \text{otherwise}
        \end{cases} \,,
\end{align}
with $u_{L1} = 1$, $u_{R1} = u_{L2} = 0.1$ and $u_{R2}=-0.5$, and the exact solution is prescribed as
\begin{align}
    \label{2D-exact_shock}
     u(x, y, t) =
        \begin{cases}
        u_{L1} & \text{if } x < -0.8 + s_{x1} t \text{ or } y < -0.8 + s_{y1} t \text{ or } x + y < -0.8 + (s_{x1} + s_{y1}) t \\
        u_{L2} & \text{if } x < 0.2 + s_{x2} t  \text{ or } y < 0.2 + s_{y2} t  \text{ or } x + y < 0.7 + (s_{x2} + s_{y2}) t \\
        u_{R2} & \text{otherwise}
        \end{cases} \,,
\end{align}
with $s_{x1}$,\,$s_{x2}$,\,$s_{y1}$ and $s_{y2}$ being the shock speeds
\begin{align}
    s_{x1} = \frac{f(u_{L1}) - f(u_{R1})}{u_{L1} - u_{R1}} \,, \hspace{0.5cm} s_{x2} = \frac{f(u_{L2}) - f(u_{R2})}{u_{L2} - u_{R2}} \,, \nonumber\\
    s_{y1} = \frac{g(u_{L1}) - g(u_{R1})}{u_{L1} - u_{R1}} \,, \hspace{0.5cm} s_{y2} = \frac{g(u_{L2}) - g(u_{R2})}{u_{L2} - u_{R2}}\,. 
\end{align} 

\begin{figure}[h]
    \begin{center}
    \includegraphics[width=\textwidth]{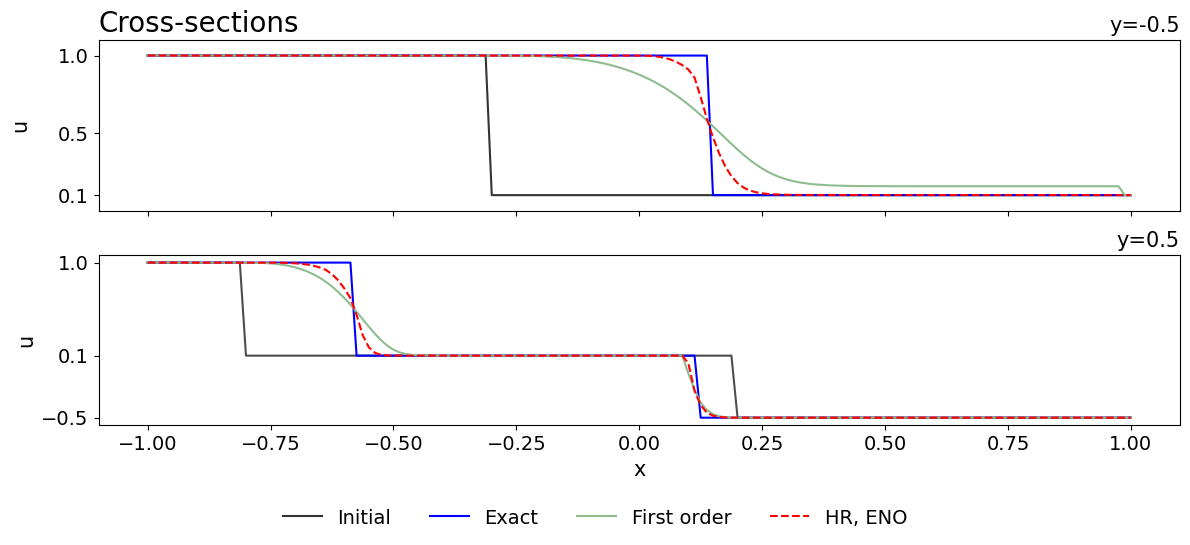}
    \end{center}
    \caption{Cross-sectional profiles at $y=-0.5$ and $y=0.5$. The plots compare the initial condition, exact solution, first order method, and the high-resolution scheme with ENO approximations, demonstrating the differences in accuracy and resolution for the example \ref{example-burgers-shocks}.}
    \label{FIG:2D-shock_crosssection}
\end{figure} 

This time, we choose $N=M/40$ number of time steps for $M = 80, 160, 320, 640$ (corresponding to 2, 4, 8, 16 time steps), with the value $C_{max} = 8$, to demonstrate the schemes also with higher Courant number. The ENO and WENO  schemes were used to compute the solution together with the first order accurate scheme, to demonstrate the difference in the errors obtained, see Table \ref{TAB:secondexample_burgers_shocks}. 
\rv{Note that, as always, these results are obtained after four sweeps \eqref{sweeps}.
}

For this case, the maximum and minimum values are not included in Table \ref{TAB:secondexample_burgers_shocks}, as no oscillations were detected and the solution consistently remains within the interval $[-0.5, 1]$.

\rv{
The residual norms, $||r||_1$ and $||r||_{\infty}$, were tracked throughout the computations and found to decrease as the mesh is refined and additional iterations are applied. Again, due to a larger complexity of this example, the rate of convergence is slower, especially for the stringent  maximum norm $||r||_{\infty}$, but it is always decreasing. 
}

These findings may be directly compared to the results for the first-order accurate scheme provided in Table \ref{TAB:secondexample_burgers_shocks}.

\begin{table}[h]
    \begin{center}
    \begin{tabular}{c c | c c c c }
    \multicolumn{2}{c}{} & \multicolumn{4}{c}{$1^{st}$ order, 4GS (8GS) [12GS]} \\
    \hline
    $M $& $N$ & E     & EOC  & $||r||_1$         & $||r||_{\infty}$         \\ 
    \hline
    40  & 1  & 0.4092 & -    & $2.0 \cdot 10^{-4}$ ($8.9 \cdot 10^{-16}$) [$2.2 \cdot 10^{-15}$] & $1.0 \cdot 10^{-2}$ ($7.9 \cdot 10^{-15}$) [$2.2 \cdot 10^{-15}$]\\
    80  & 2  & 0.2584 & 0.66 & $5.9 \cdot 10^{-5}$ ($1.5 \cdot 10^{-15}$) [$8.2 \cdot 10^{-15}$] & $9.7 \cdot 10^{-3}$ ($8.8 \cdot 10^{-15}$) [$8.2 \cdot 10^{-15}$]\\
    160 & 4  & 0.1527 & 0.75 & $5.4 \cdot 10^{-5}$ ($4.3 \cdot 10^{-8}$) [$2.4 \cdot 10^{-10}$] & $9.8 \cdot 10^{-3}$  ($1.3 \cdot 10^{-4}$) [$2.4 \cdot 10^{-6}$]\\
    320 & 8  & 0.0841 & 0.85 & $2.1 \cdot 10^{-5}$ ($2.6 \cdot 10^{-8}$) [$8.8 \cdot 10^{-10}$] & $9.8 \cdot 10^{-3}$  ($3.2 \cdot 10^{-4}$) [$8.8 \cdot 10^{-6}$]\\
    \hline \\
    \multicolumn{2}{c}{} & \multicolumn{4}{c}{ENO, 4GS (8GS) [12GS]} \\
    \hline
    $M $& $N$& E     & EOC  & $||r||_1$         & $||r||_{\infty}$\\ 
    \hline
    40  & 1  & 0.2233 & -    & $5.3 \cdot 10^{-5}$ ($1.5 \cdot 10^{-9}$) [$2.4 \cdot 10^{-13}$] & $6.1 \cdot 10^{-3}$ ($3.2 \cdot 10^{-7}$) [$4.6 \cdot 10^{-11}$]\\
    80  & 2  & 0.1231 & 0.85 & $9.9 \cdot 10^{-5}$ ($1.5 \cdot 10^{-7}$) [$1.8 \cdot 10^{-9}$] & $1.4 \cdot 10^{-2}$ ($1.1 \cdot 10^{-4}$) [$2.5 \cdot 10^{-6}$]\\
    160 & 4  & 0.0649 & 0.92 & $7.3 \cdot 10^{-5}$ ($2.3 \cdot 10^{-8}$) [$2.3 \cdot 10^{-10}$] & $2.1 \cdot 10^{-2}$ ($7.1 \cdot 10^{-5}$) [$1.2 \cdot 10^{-6}$]\\
    320 & 8  & 0.0335 & 0.95 & $2.5 \cdot 10^{-5}$ ($3.3 \cdot 10^{-15}$) [$3.2 \cdot 10^{-15}$] & $1.1 \cdot 10^{-2}$ ($1.7 \cdot 10^{-12}$) [$1.2 \cdot 10^{-14}$]\\
    \hline \\
    \multicolumn{2}{c}{} & \multicolumn{4}{c}{WENO, 4GS (8GS) [12GS]} \\
    \hline
    $M $& $N$& E     & EOC  & $||r||_1$         & $||r||_{\infty}$  \\ 
    \hline
    40  & 1  & 0.2389 & -    & $1.1 \cdot 10^{-4}$ ($4.0 \cdot 10^{-7}$) [$1.9 \cdot 10^{-9}$] & $1.0 \cdot 10^{-2}$ ($1.2 \cdot 10^{-4}$) [$1.9 \cdot 10^{-6}$]\\
    80  & 2  & 0.1425 & 0.74 & $9.1 \cdot 10^{-5}$ ($5.6 \cdot 10^{-7}$) [$7.7 \cdot 10^{-8}$] & $1.2 \cdot 10^{-2}$ ($4.1 \cdot 10^{-4}$) [$7.7 \cdot 10^{-6}$]\\
    160 & 4  & 0.0776 & 0.87 & $4.6 \cdot 10^{-5}$ ($1.3 \cdot 10^{-7}$) [$7.3 \cdot 10^{-9}$] & $2.7 \cdot 10^{-2}$ ($3.8 \cdot 10^{-4}$) [$7.3 \cdot 10^{-6}$]\\
    320 & 8  & 0.0399 & 0.95 & $2.8 \cdot 10^{-5}$ ($7.9 \cdot 10^{-8}$) [$1.1 \cdot 10^{-9}$] & $4.4 \cdot 10^{-2}$ ($5.6 \cdot 10^{-4}$) [$1.3 \cdot 10^{-5}$]\\
    \hline 
    \end{tabular}
    \end{center} 
    \caption{The numerical errors and the EOCs for the first and second (ENO, WENO) order accurate schemes for the nonlinear problem with the initial condition \eqref{2D-exact_shock}, $T=0.4$, $C_{max} = 8$. }
    \label{TAB:secondexample_burgers_shocks}
\end{table}

\clearpage
\subsubsection{\rv{Fourth example}}

As the last example, we show results of a combination of Burgers' equation with the velocity field \eqref{velocityfield}
\begin{equation}
    \label{Burgers2D}
    \partial_{t} u +  \partial_x \Big( v \frac{u^2}{2}\Big) +  \partial_y \Big(w \frac{u^2}{2}\Big) = 0 \,.
\end{equation}

For this case, we employ the initial condition defined in \eqref{InitialSpec}, which consists of four distinct shapes in the four quadrants of the computational domain $x,y \in (-1.25,1.25)$. We simulate a segment of the rotation with $T=1/8$. This case necessitates the use of the high-resolution scheme with used number of finite volumes $M=20, 40, 80, 160$ and the number of time steps $N=M/10$ being 2, 4, 8, 16. Given the maximum function value being 1, the maximum Courant number is bounded as $C_{max}^x = C_{max}^y = 3.92$.


Since the exact solution is unknown, we omit error and EOC analysis. Instead, the Table \ref{TAB:fourthexample_burgerswithrotation} lists the minimum and maximum values obtained numerically with the first order accurate scheme and the high-resolution schemes with ENO and WENO approximations. Computations are performed using four and eight Gauss-Seidel iterations, with the latter marginally improving the residuals toward numerical zero. 

\begin{table}[h]
    \begin{center}
    \begin{tabular}{c c | c c c c}
    \multicolumn{2}{c}{} & \multicolumn{4}{c}{$1^{st}$ order, 4GS (8GS)} \\
    \hline
    $M $& $N$ & min     & max            & $||r||_1$         & $||r||_{\infty}$         \\ 
    \hline
    20  & 2   & 0.0 & 0.71 & $2.0 \cdot 10^{-8}$  ($4.5 \cdot 10^{-17}$)& $1.2 \cdot 10^{-6}$  ($2.2 \cdot 10^{-16}$)\\
    40  & 4   & 0.0 & 0.77 & $2.2 \cdot 10^{-10}$ ($4.3 \cdot 10^{-17}$)& $3.4 \cdot 10^{-8}$  ($2.4 \cdot 10^{-16}$)\\
    80  & 8   & 0.0 & 0.86 & $3.8 \cdot 10^{-14}$ ($4.7 \cdot 10^{-17}$)& $1.3 \cdot 10^{-12}$ ($6.1 \cdot 10^{-16}$)\\
    160 & 16  & 0.0 & 0.93 & $3.4 \cdot 10^{-14}$ ($5.1 \cdot 10^{-17}$)& $1.8 \cdot 10^{-12}$ ($7.7 \cdot 10^{-16}$)\\
    \hline \\
    \multicolumn{2}{c}{} & \multicolumn{4}{c}{ENO, 4GS (8GS)} \\
    \hline
    $M $& $N$& min                   & max  & $||r||_1$         & $||r||_{\infty}$ \\ 
    \hline
    20  & 2  & $-9.8 \cdot 10^{-118}$& 0.82 & $1.2 \cdot 10^{-9}$  ($3.4 \cdot 10^{-17}$)& $3.2 \cdot 10^{-8}$ ($1.6 \cdot 10^{-16}$)  \\
    40  & 4  & $-4.8 \cdot 10^{-45}$ & 0.92 & $1.4 \cdot 10^{-14}$ ($3.9 \cdot 10^{-17}$)& $3.6 \cdot 10^{-12}$ ($2.7 \cdot 10^{-16}$)  \\
    80  & 8  & $-1.4 \cdot 10^{-59}$ & 0.97 & $4.3 \cdot 10^{-17}$ ($4.3 \cdot 10^{-17}$)& $4.4 \cdot 10^{-16}$ ($3.3 \cdot 10^{-16}$)  \\
    160 & 16 & $-7.7 \cdot 10^{-30}$ & 0.99 & $4.3 \cdot 10^{-17}$ ($4.4 \cdot 10^{-17}$)& $4.4 \cdot 10^{-16}$ ($5.5 \cdot 10^{-16}$)  \\
    \hline \\
    \multicolumn{2}{c}{} & \multicolumn{4}{c}{WENO, 4GS (8GS)} \\
    \hline
    $M $& $N$& min  & max & $||r||_1$         & $||r||_{\infty}$  \\ 
    \hline
    20  & 2  & $-1.7 \cdot 10^{-21}$ & 0.83 & $1.3 \cdot 10^{-9}$ ($3.7 \cdot 10^{-17}$)& $2.7 \cdot 10^{-8}$ ($1.6 \cdot 10^{-16}$)  \\
    40  & 4  & $-9.1 \cdot 10^{-21}$ & 0.93 & $1.7 \cdot 10^{-14}$ ($4.0 \cdot 10^{-17}$)& $4.1 \cdot 10^{-12}$ ($3.3 \cdot 10^{-16}$)  \\
    80  & 8  & $-2.3 \cdot 10^{-18}$ & 0.98 & $4.1 \cdot 10^{-17}$ ($4.1 \cdot 10^{-17}$)& $3.8 \cdot 10^{-16}$ ($3.8 \cdot 10^{-16}$)  \\
    160 & 16 & $-1.6 \cdot 10^{-18}$ & 0.99 & $4.0 \cdot 10^{-17}$ ($4.0 \cdot 10^{-17}$)& $5.5 \cdot 10^{-16}$ ($5.5 \cdot 10^{-16}$)  \\
    \hline 
    \end{tabular}
    \end{center}
    \caption{The minimum and maximum values of the first order scheme and the high-resolution scheme with ENO and WENO approximations for the combination of Burgers' equation with the rotation of four shapes \eqref{InitialSpec} with $C_{max}^x = C_{max}^y = 3.92$. The norms of residuals are shown for four (4GS) and eight (8GS) Gauss-Seidel iterations. }
    \label{TAB:fourthexample_burgerswithrotation}
\end{table}

\begin{figure}[h]
    \begin{center}
    \includegraphics[width=\textwidth]{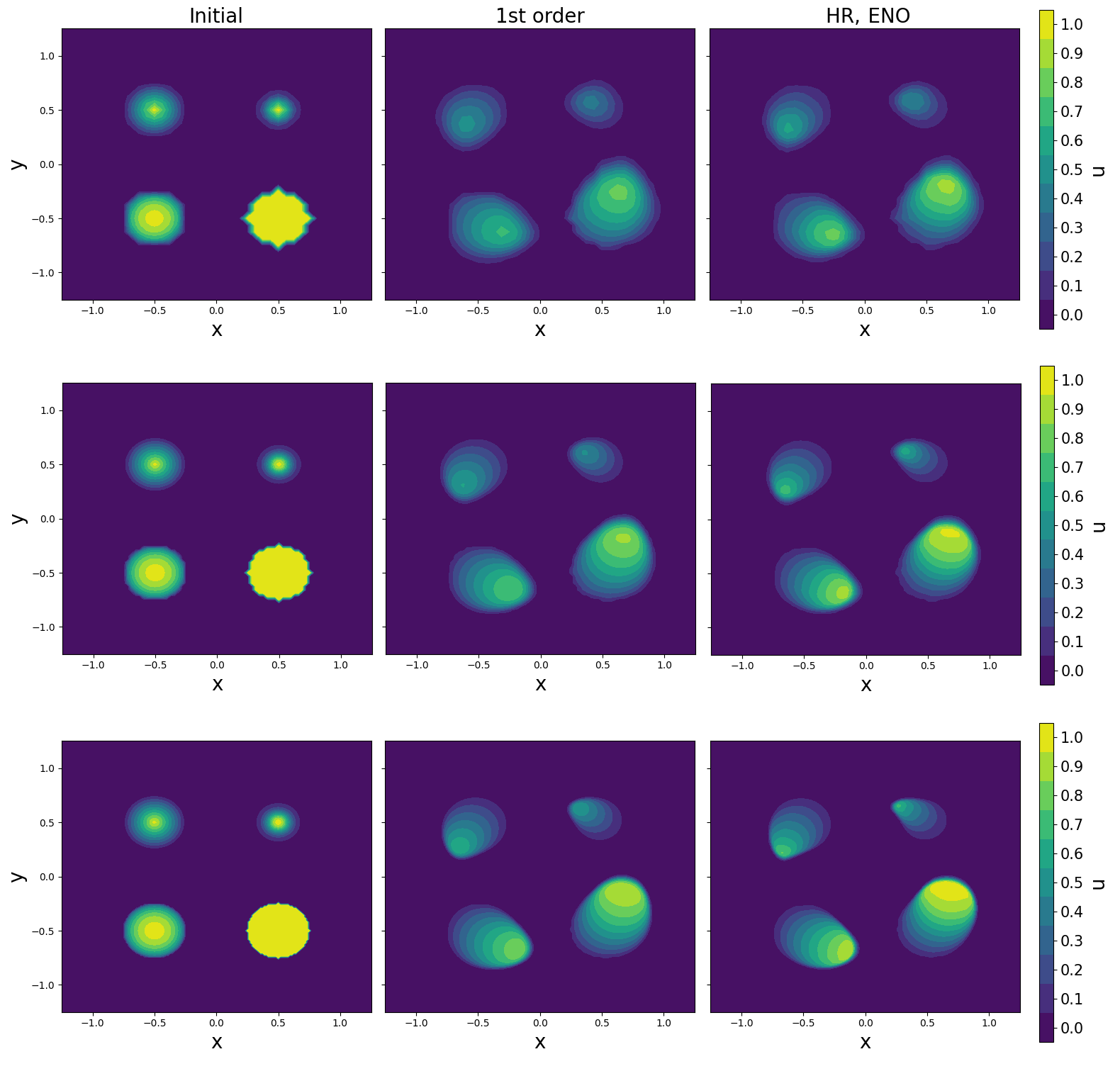}
    \end{center}
    \caption{The comparison of the contours of the numerical solution obtained with the first order (the middle column) and the high-resolution with ENO approximation (the right column) for the different choices of $M = 40$ (top row), $M = 80$ (middle row), $M = 160$ (bottom row) for the final time $T=1/8$. The initial condition is shown in the left column. }
\label{FIG:fourthexample_burgerswithrotation}
\end{figure}

To compensate for the lack of an exact solution, we provide visualizations of the progression of numerical solutions under mesh refinement $(M=40 (N=4), M=80 (N=8), 160 (N=16))$. The Figure \ref{FIG:fourthexample_burgerswithrotation} displays contour plots that compare the results from the first order scheme and the high-resolution scheme (ENO). These figures highlight significant qualitative improvements when using the high-resolution method and the finer meshes, particularly in preserving sharp features and reducing numerical dissipation.

\clearpage
\section{Conclusion}
\label{sec-conclude}

In this paper, we derived the high-resolution compact implicit finite volume numerical schemes in two dimensions for solving hyperbolic problems in the form of scalar conservation laws. Our approach utilizes the compactness of the stencil for computational efficiency, as the solution of resulting algebraic systems can be efficiently obtained through the fast sweeping method combined with nonlinear Gauss-Seidel iterations.

The parametric form of the second order scheme inherently supports the use of Essentially Non-Oscillatory (ENO) and Weighted Essentially Non-Oscillatory (WENO) approximations for handling discontinuous solutions. 
\rv{
We proposed the limiting in time that is used with such unmodified (W)ENO approximations and that is provably delivering oscillations free numerical solutions for linear advection with variable velocity if some standard properties of the space reconstruction are fulfilled.
}

We validated these numerical schemes through numerical experiments on selected examples of linear and nonlinear scalar conservation laws demonstrating their properties in various scenarios.



\printbibliography
\end{document}